\documentclass[10pt]{book}
\usepackage[sectionbib]{natbib}
\usepackage{array,epsfig,fancyheadings,rotating}

\usepackage{amsmath}
\usepackage{amssymb}
\usepackage{amsfonts}
\usepackage{multirow}
\usepackage{amsthm}

\newtheorem{theorem}{Theorem}
\newtheorem{lemma}{Lemma}
\newtheorem{corollary}{Corollary}

\theoremstyle{definition}

\newtheorem{example}{Example}
\newtheorem{remark}{Remark}
\usepackage{amsbsy}
\usepackage{setspace}
\newtheorem{assumption}{Assumption}

\renewcommand{\thelemma}{\arabic{lemma}}

\begin{document}

	\markboth{\hfill{\footnotesize\rm JIE YANG, LIPING TONG AND ABHYUDAY MANDAL} \hfill}
	{\hfill {\footnotesize\rm D-OPTIMAL DESIGNS WITH ORDERED CATEGORICAL DATA} \hfill}
	
	\renewcommand{\thefootnote}{}
	$\ $\par


\thispagestyle{empty}

\fontsize{12}{14pt plus.8pt minus .6pt}\selectfont \vspace{0.8pc}
\centerline{\large\bf D-OPTIMAL DESIGNS WITH ORDERED CATEGORICAL DATA}
\vspace{.4cm} \centerline{Jie Yang$^{1}$, Liping Tong$^{2}$ and Abhyuday Mandal$^{3}$} \vspace{.4cm} \centerline{\it
$^1$University of Illinois at Chicago, $^2$Advocate Health Care and $^3$University of Georgia} \vspace{.55cm} \fontsize{9}{11.5pt plus.8pt minus
.6pt}\selectfont

\begin{quotation}
\noindent {\it Abstract:}
Cumulative link models have been widely used for ordered categorical responses. Uniform allocation of experimental units is commonly used in practice, but often suffers from a lack of efficiency. We consider D-optimal designs with ordered categorical responses and cumulative link models. For a predetermined set of design points, we derive the necessary and sufficient conditions for an allocation to be locally D-optimal and develop efficient algorithms for obtaining approximate and exact designs. We prove that the number of support points in a minimally supported design only depends on the number of predictors, which can be much less than the number of parameters in the model. We show that a D-optimal minimally supported allocation in this case is usually not uniform on its support points. In addition, we provide EW D-optimal designs as a highly efficient surrogate to Bayesian D-optimal designs. Both of them can be much more robust than uniform designs.
\par

\vspace{9pt}
\par
\end{quotation}\par

\def\thefigure{\arabic{figure}}
\def\thetable{\arabic{table}}

\renewcommand{\theequation}{\thechapter.\arabic{equation}}

\fontsize{12}{14pt plus.8pt minus .6pt}\selectfont

\setcounter{chapter}{1}
\setcounter{section}{1}
\setcounter{equation}{0} 
\noindent  {\bf 1. Introduction}

In this paper we determine optimal and efficient designs for factorial experiments with qualitative factors and ordered categorical responses, or simply ordinal data. Design of experiment with multinomial response, and ordered categories in particular, is becoming increasingly popular in a rich variety of scientific disciplines, especially when human evaluations are involved (\cite{christensen2013}). Examples include a wine bitterness study (\cite{randall1989}), potato pathogen experiments (\cite{omer2000}), a radish seedling's damping-off study (\cite{krause2001}), a polysilicon deposition study (\cite{wu2008}), beef cattle research (\cite{osterstock2010}), and a toxicity study (\cite{agresti2013}).

This research is motivated by an {\it odor removal study} conducted by the textile engineers at the University of Georgia.
The scientists studied the manufacture of bio-plastics containing odorous volatiles, that need to be removed before commercialization. For that purpose, a $2^2$ factorial experiment was conducted using algae and synthetic plastic resin blends. The factors were {\tt types of algae} ($x_1$: raffinated or solvent extracted algae ($-$), catfish pond algae ($+$)) and {\tt synthetic resins} ($x_2$: polyethylene ($-$), polypropylene ($+$)). The response $Y$ had three ordered categories: serious odor ($j=1$), medium odor ($j=2$), and almost no odor ($j=3$). Following traditional factorial design theory, a pilot study with equal numbers ($10$ in this case) of replicates at each experimental setting was conducted, a {\it uniform design}.  The results are summarized in Table~\ref{tab:odor}, where $y_{ij}$ represents the number of responses falling into the $j$th category under the $i$th experimental setting.
As demonstrated later (Section~4), the best design identified by our research could improve the efficiency by 25\% with only three experimental settings involved.

\lhead[\footnotesize\thepage\fancyplain{}\leftmark]{}\rhead[]{\fancyplain{}\rightmark\footnotesize\thepage}

\begin{table}[h!]\caption{Pilot Study of Odor Removal Study}\label{tab:odor}
\begin{center}
\begin{tabular}{|c|cc|ccc|}
\hline
Experimental & Factor & level  & \multicolumn{3}{c|}{Summarized responses ($Y$, odor)} \\
setting & Algae & Resin & Serious & Medium & No odor  \\
$i$ & $x_1$ & $x_2$ & $y_{i1}$ & $y_{i2}$ & $y_{i3}$  \\
\hline
$1$ & $+$ & $+$ &  2 & 6 &  2   \\
$2$ & $+$ & $-$ &  7 & 2 &  1   \\
$3$ & $-$ & $+$ &  0 & 0 & 10   \\
$4$ & $-$ & $-$ &  0 & 2 &  8   \\
\hline
\end{tabular}
\end{center}
\end{table}

For such kind of ordinal response $Y$ with $J$ categories and $d$ predictors ${\mathbf x}=(x_1, \ldots, x_d)^T$, the most popular model in practice was first the {\it proportional odds model} (also known as {\it cumulative logit model}, see \cite{agresti2005} for a detailed review). \cite{pmcc1980} extended it to the {\it cumulative link model} (also known as {\it ordinal regression model})
\begin{equation}\label{clm}
g\left(P(Y\leq j\mid {\mathbf x})\right)=\theta_j - \boldsymbol\beta^T{\mathbf x}, \ j=1, \ldots, J-1
\end{equation}
where $g$ is a general link function, with the proportional odds model as a special case when $g$ is the logit link. Examples include the complementary log-log link for the polysilicon deposition study (see Example~\ref{industrialexample}) and the cauchit link for the toxicity study (see Example~\ref{d=1J=3m=5example}). We adopt the cumulative link model (\ref{clm}).

When there are only two categories ($J=2$), the cumulative link model (\ref{clm}) is essentially a generalized linear model for binary data (McCullagh and Nelder (1989); Dobson and Barnett (2008)).
For optimal designs under generalized linear models, there is a growing body of literature (see \cite{khuri2006}, \cite{atkinson2007}, \cite{stufken2012}, and references therein). In this case, it is known that the minimum number of experimental settings required by a nondegenerate Fisher information matrix is $d+1$, which equals the number of parameters (\cite{fedorov1972, ym2014}). A design with the least number of experimental settings, known as a {\it minimally supported design}, is of practical significance with a specified regression model due to the cost of changing settings. It is also known that the experimental units should be uniformly assigned when a minimally supported design is adopted for binary response, or under a univariate generalized linear model (\cite{ym2014}).

When $J\geq 3$, the cumulative link model is a special case of the multivariate generalized linear model (\cite{pmcc1980}). The relevant results in the optimal design literature are meagre and restricted to the logit link function (\cite{atkinson1999, perevozskaya2003}).
Here we obtain theoretical results and efficient algorithms for general link functions and reveal that the optimal designs with $J\geq 3$ are quite different from the cases with $J=2$. We prove that the minimum number of experimental settings is still $d+1$, but strictly less than the number of parameters $d+J-1$ (Theorems~\ref{mindesigntheorem} and \ref{mindesigntheorem2}). This counter-intuitive result is due to the multinomial-type responses: from a single experimental setup, the summarized responses have $J-1$ degrees of freedom, requiring fewer distinct experimental settings in a minimally supported design. For the same reason, the allocation of replicates in a minimally supported design is usually not uniform (Section~5), which differs from the traditional factorial design theory.

As with generalized linear models, the information matrix under cumulative link models depends on unknown parameters. Different approaches have been proposed to solve the dependence of optimal designs on unknown parameters, including local optimality (\cite{chernoff1953}), Bayesian approach (\cite{chaloner1995}), a maximin approach
(Pronzato and Walter (1988); Imhof (2001)), 
and a sequential procedure (\cite{ford1989}).
As pointed out by \cite{ford1992}, locally optimal designs are not only important when good initial parameters are available from previous experiments, but can also be a benchmark for designs chosen to satisfy experimental constraints. We mainly focus on locally optimal designs. For situations where local values of the parameters are difficult to obtain, but the experimenter has an idea of the range of parameters with or without a prior distribution, we recommend EW optimal designs, where the Fisher information matrix is replaced by its expected values (\cite{atkinson2007,ymm2013}). We compare Bayesian D-optimal designs (\cite{chaloner1995}) with EW D-optimal designs for ordinal data. As a surrogate for Bayesian designs, an EW design is much easier to find and retains high efficiency with respect to Bayesian criterion (Section~6).

Among various optimal design criteria, D-optimality, which maximizes the determinant of Fisher information matrix, is the most frequently used (\cite{atkinson1999}) and often performs well according to other criteria (\cite{atkinson2007}). We study D-optimal designs.

In the design literature, one type of experiment deals with quantitative or continuous factors only. Such a design problem includes the identification of a set of design points $\{{\bf x}_i\}_{i=1,\ldots, m}$ and the corresponding weights $\{ p_i\}_{i=1, \ldots, m}$ (see, for example, \cite{atkinson2007} and \cite{stufken2012}). Numerical algorithms are typically used for cases with two or more factors (see, for example, \cite{woods2006}). Another type of experiment employs qualitative or discrete factors, where the set of design points $\{{\bf x}_i\}_{i=1,\ldots, m}$ is predetermined and only the weights $\{p_i\}_{i=1, \ldots, m}$ are to be optimized (see, for example, \cite{ym2014}). One can pick grid points of continuous factors and turn the first kind of problem into the second. \citeauthor{tong2014}~(2014, Section~5) also bridged the gap between the two types of problems in a way that results involving discrete factors can be applied to the cases with continuous factors. We concentrate on the second kind of design problems and assume that $\{{\bf x}_i\}_{i=1,\ldots, m}$ are given and fixed.

This paper is organized as follows.  In Section~2, we describe the preliminary setup and obtain the Fisher information matrix for the cumulative link model with a general link, generalizing \cite{perevozskaya2003}. We also identify a necessary and sufficient condition for the Fisher information matrix to be positive definite, which determines the minimum number of experimental settings required. In Sections~3 and 4, we provide theoretical results and numerical algorithms for searching locally D-optimal approximate or exact designs. In Section~5, we identify analytic D-optimal designs for special cases to illustrate that a D-optimal minimally supported design is usually not uniform on its support points. In Section~6, we illustrate by examples that the EW D-optimal design can be highly efficient with respect to Bayesian D-optimality. We make concluding remarks in Section~7 and relegate additional proofs and results to the supplementary materials.


\setcounter{chapter}{2}
\setcounter{section}{2}
\setcounter{equation}{0} 

\smallskip
\noindent {\bf 2. Fisher Information Matrix and Its Determinant}

Suppose there are $m$ ($m\geq 2$) predetermined experimental settings. For the $i$th experimental setting with corresponding predictors ${\mathbf x}_i = (x_{i1},$ $\ldots,$ $x_{id})^T$ $\in$ $\mathbb{R}^d$ ($d\geq 1$), there are $n_i$ experimental units assigned to it. Among the $n_i$ experimental units, the $k$th one generates a response $V_{ik}$ which belongs to one of $J$ ($J \geq 2$) ordered categories. As shown in Example~\ref{profitd3abexample}, the dimension $d$ of the predictors can be significantly larger than the number of factors considered in the experiment, which allows more flexible models.


\noindent   {\bf 2.1 General setup}

In many applications, $V_{i1}, \ldots, V_{in_i}$ are regarded as i.i.d.~discrete random variables. Let $\pi_{ij} = P(V_{ik}=j)$, where $i=1, \ldots, m$; $j=1, \ldots, J$; and $k=1, \ldots, n_i$~. Let $Y_{ij} = \#\{k\mid V_{ik}=j\}$ be the number of $V_{ik}$'s falling into the $j$th category. Then $(Y_{i1}, \ldots, Y_{iJ}) \sim {\rm Multinomial} (n_i; \pi_{i1}, \ldots, \pi_{iJ})$~.
\begin{assumption}\label{pi>0}
$0<\pi_{ij}<1$, $i=1, \ldots, m$; $j=1, \ldots, J$.
\end{assumption}

Let $\gamma_{ij}=P(V_{ik} \leq j) = \pi_{i1}+\cdots + \pi_{ij}\ ,\>\> j=1, \ldots, J.$ Based on Assumption~\ref{pi>0}, $0<\gamma_{i1}<\gamma_{i2}<\cdots<\gamma_{i,J-1}<\gamma_{iJ}=1$ for each $i=1, \ldots, m$. Consider independent multinomial observations $(Y_{i1}, \ldots, Y_{iJ}), i=1,$ $\ldots,$ $m$ with corresponding predictors ${\mathbf x}_1, \ldots, {\mathbf x}_m$~. Under a {\it cumulative link model} or {\it ordinal regression model} (\cite{pmcc1980,agresti2013,christensen2013}), there exists a link function $g$ and parameters of interest $\theta_1, \ldots, \theta_{J-1}, \boldsymbol{\beta} = (\beta_1, \ldots, \beta_d)^T$, such that
\begin{eqnarray}\label{eq:gamma}
g(\gamma_{ij}) = \theta_j - {\mathbf x}_i^T\boldsymbol{\beta}, \quad j=1, \ldots, J-1.
\end{eqnarray}
This leads to $m(J-1)$ equations in $d+J-1$ parameters $(\beta_1, \ldots,$ $\beta_d,$ $\theta_1, \ldots,$ $\theta_{J-1})$.

\begin{assumption}\label{glink}
The link $g$ is differentiable and its derivative $g'>0$.
\end{assumption}

Assumption~\ref{glink} is satisfied for commonly used link functions including {\tt logit} ($\log (\gamma/(1-\gamma))$, {\tt probit} ($\Phi^{-1}(\gamma)$), {\tt log-log} ($-\log(-\log(\gamma))$), {\tt complementary log-log} ($\log(-\log(1-\gamma))$), and {\tt cauchit} ($\tan(\pi(\gamma-1/2))$) (\cite{pmcc1989, christensen2013}). Some relevant formulas of these link functions are provided in the supplementary materials (Section~\ref{S1section}). According to Assumption~\ref{glink}, $g$ is strictly increasing, and then $\theta_1 < \theta_2 < \cdots <\theta_{J-1}$~.

\begin{example}\label{logitd2J3example}
{\rm
Consider the logit link $g(\gamma)=\log(\gamma/(1-\gamma))$ with two predictors and three ordered categories. Model~(\ref{eq:gamma}) consists of $2m$ equations
$g(\gamma_{ij})=\theta_j - x_{i1}\beta_1 - x_{i2}\beta_2, \>\>i=1,\ldots, m; \ j=1,2$
and parameters $(\beta_1, \beta_2, \theta_1, \theta_2)$. Under Assumptions~\ref{pi>0} and \ref{glink}, $\theta_1 < \theta_2$~.
\hfill{$\Box$}}\end{example}

\begin{example}\label{profitd3abexample}
{\rm
Suppose the model consists of three covariates $x_1, x_2, x_3$ and a few second-order predictors,
$g(\gamma_{ij})=\theta_j - x_{i1}\beta_1 - x_{i2}\beta_2 - x_{i3}\beta_3 - x_{i1}x_{i2}\beta_{12} - x_{i1}^2\beta_{11} - x_{i2}^2\beta_{22}$,
where $i=1,\ldots, m; j=1,\ldots, J-1$. Then the number of predictors is $d=6$.
\hfill{$\Box$}}\end{example}

Under the cumulative link model~(\ref{eq:gamma}), the log-likelihood function (up to a constant) is
$l(\beta_1, \ldots, \beta_d, \theta_1, \ldots,\theta_{J-1}) = \sum_{i=1}^m \sum_{j=1}^J Y_{ij}\log (\pi_{ij})$,
where $\pi_{ij} = \gamma_{ij} - \gamma_{i,j-1}$ with $\gamma_{ij} = g^{-1}(\theta_j - {\mathbf x}_i^T \boldsymbol{\beta})$ for $j=1,\ldots, J-1$ and $\gamma_{i0}=0$, $\gamma_{iJ}=1$.

\cite{perevozskaya2003} obtained a detailed form of the Fisher information matrix for logit link and one predictor. Our result is for general link and $d$ predictors; its proof is relegated to the supplementary materials (Section~\ref{S3section}).

\begin{theorem}\label{sumFtheorem}
Under Assumptions \ref{pi>0} and \ref{glink}, the Fisher information matrix can be written as
\begin{equation}\label{Fsum}
{\mathbf F}=\sum_{i=1}^m n_i {\mathbf A}_i
\end{equation}
where ${\mathbf A}_i$ is the $(d+J-1)\times (d+J-1)$ matrix
\[
\left(
\begin{array}{cc}
{\mathbf A}_{i1} & {\mathbf A}_{i2} \\
{\mathbf A}_{i2}^T & {\mathbf A}_{i3}\\
\end{array}
\right)
=
\left(
\begin{array}{cc}
(e_i x_{is}x_{it})_{s=1,\ldots d; t=1,\ldots,d} & (-x_{is}c_{it})_{s=1,\ldots,d; t=1,\ldots,J-1} \\
(-c_{is}x_{it})_{s=1,\ldots,J-1; t=1,\ldots,d} & {\mathbf A}_{i3}\\
\end{array}
\right)
\]
and ${\mathbf A}_{i3}$ is the $(J-1)\times (J-1)$ symmetric tri-diagonal matrix with diagonal entries $u_{i1}, \ldots,$ $u_{i, J-1}$, and off-diagonal entries $-b_{i2}, \ldots, -b_{i,J-1}$ when $J\geq 3$, where
$e_i = \sum_{j=1}^J \pi_{ij}^{-1} (g_{ij}-g_{i,j-1})^2 > 0$ with $g_{ij}=(g^{-1})'(\theta_j-{\mathbf x}_i^T\boldsymbol{\beta})>0$ for $j=1, \ldots, J-1$ and $g_{i0}=g_{iJ}=0$; $c_{it} = g_{it}[\pi_{it}^{-1} (g_{it}-g_{i,t-1}) - \pi_{i,t+1}^{-1} (g_{i,t+1}-g_{it})]$; $u_{it} = g_{it}^2(\pi_{it}^{-1}+\pi_{i,t+1}^{-1}) > 0$; and $b_{it} = g_{i,t-1}g_{it}\pi_{it}^{-1} > 0$.
${\mathbf A}_{i3}$ contains only one entry $u_{i1}$ when $J=2$.
\end{theorem}

As the Fisher information matrix, ${\mathbf F}$ is always positive semi-definite, $|{\mathbf F}| \geq 0$ (\cite{fedorov1972}). As a special case, ${\mathbf A}_i$ is the Fisher information at the experimental setting ${\mathbf x}_i$ (also known as a {\it design point} or {\it support point}) and thus is positive semi-definite.


\medskip
\noindent   {\bf 2.2 Determinant of Fisher information matrix}

Among different criteria for optimal designs, D-criterion looks for the allocation maximizing $|{\mathbf F}|$, the determinant of ${\mathbf F}$. Here, a D-optimal design with $m$ predetermined design points ${\mathbf x}_1, \ldots, {\mathbf x}_m$ could either be an integer-valued allocation $(n_1,$ $n_2,$ $\ldots,$ $n_m)$ maximizing $|{\mathbf F}|$ with fixed $n = \sum_{i=1}^m n_i > 0$, known as an {\it exact design}; or a real-valued allocation $(p_1, p_2, \ldots, p_m)$ maximizing $|n^{-1}{\mathbf F}|$ with $p_i = n_i/n \geq 0$ and $\sum_{i=1}^m p_i=1$, known as an {\it approximate design}.

\begin{theorem}\label{polynomialtheorem}
The determinant of the Fisher information matrix,
$$|{\mathbf F}|=\sum_{\alpha_1 + \cdots + \alpha_m=d+J-1} c_{\alpha_1, \ldots, \alpha_m}\cdot n_1^{\alpha_1}\cdots n_m^{\alpha_m}\ ,$$
is an order-$(d+J-1)$ homogeneous polynomial of $(n_1, \ldots, n_m)$ and
\begin{equation}\label{Atau}
c_{\alpha_1, \ldots, \alpha_m} = \sum_{\tau \in (\alpha_1, \ldots, \alpha_m)}|{\mathbf A}_\tau|\ .
\end{equation}
\end{theorem}

The proof of Theorem~\ref{polynomialtheorem} is relegated to the supplementary materials (Section~\ref{S3section}). Given a map $\tau : \{1, 2, \ldots, d+J-1\} \rightarrow \{1, \ldots, m\}$, ${\mathbf A}_\tau$ in (\ref{Atau}) is a $(d+J-1)\times (d+J-1)$ matrix whose $k$th row is the same as the $k$th row of ${\mathbf A}_{\tau(k)}$~, $k=1, \ldots, d+J-1$. We take $\tau\in (\alpha_1, \ldots, \alpha_m)$ where $\alpha_i = \#\{j : \tau(j)=i\}$ for each $i=1, \ldots, m$.

In order to obtain analytic properties of $|{\mathbf F}|$, we need some lemmas. The first of them covers Lemma~1 in \cite{perevozskaya2003} as a special case:

\begin{lemma}\label{ranklemma}
${\rm Rank}({\mathbf A}_i)={\rm Rank}({\mathbf A}_{i3})=J-1$. Furthermore,
${\mathbf A}_{i3}$ is positive definite and
\[
|{\mathbf A}_{i3}|=\prod_{s=1}^{J-1}g_{is}^2 \cdot \prod_{t=1}^J \pi_{it}^{-1}>0
\]
where $g_{is}=(g^{-1})'(\theta_s-{\mathbf x}_i^T\boldsymbol{\beta})>0$ for $s=1, \ldots, J-1$.
\end{lemma}

\begin{example}\label{d2J3example}
{\rm
Suppose $d=2$, $J=3$, with link function $g$. According to Theorem~\ref{polynomialtheorem}, $|{\mathbf F}|$ is then an order-$4$ homogeneous polynomial of $(n_1, \ldots, n_m)$. Based on Lemma~\ref{coelemma1} and Lemma~\ref{coelemma2} in the supplementary materials (Section~\ref{S2section}), we can remove all the terms of the form $n_i^4$, $n_i^3n_j$, or $n_i^2n_j^2$ from $|{\mathbf F}|$. Therefore,
$$|{\mathbf F}|=\sum_{i=1}^m \sum_{j<k,j\neq i, k\neq i} c_{ijk}\cdot n_i^2 n_j n_k + \sum_{i<j<k<l} c_{ijkl}\cdot n_i n_j n_k n_l
$$
for some coefficients $c_{ijk}$ and $c_{ijkl}$~.
\hfill{$\Box$}}\end{example}

Based on Lemmas~\ref{coelemma1} and \ref{coelemma2}, in order to keep $c_{\alpha_1, \ldots, \alpha_m}\neq 0$, the largest possible $\alpha_i$ is $J-1$ and the fewest possible number of positive $\alpha_i$'s is $d+1$.

\begin{theorem}\label{mindesigntheorem}
$|{\mathbf F}|>0$ only if $\ m \geq d+1$.
\end{theorem}

To determine whether $d+1$ experimental settings or support points are enough to keep the Fisher information matrix positive definite, we study the leading term of $|{\mathbf F}|$ with $\max_{1\leq i\leq m} \alpha_i = J-1$. For example, $a_{i_0}=J-1$ for some $1\leq i_0 \leq m$. From Lemma~\ref{coelemma2} and $\sum_{i=1}^m \alpha_i = d+J-1$, to have $c_{\alpha_1, \ldots, \alpha_m} \neq 0$, there must exist $1\leq i_1 < i_2 < \cdots < i_d \leq m$ which are different from $i_0$, such that, $\alpha_{i_1} = \cdots = \alpha_{i_d} = 1$. A lemma provides an explicit formula for such a coefficient $c_{\alpha_1, \ldots, \alpha_m}$:
\begin{lemma}\label{coelemma3}
Suppose $\alpha_{i_0}=J-1$ and $\alpha_{i_1} = \cdots = \alpha_{i_d} = 1$. Then
$$c_{\alpha_1, \ldots, \alpha_m} = \prod_{s=1}^d e_{i_s}\cdot |{\mathbf A}_{i_03}|\cdot |{\mathbf X}_{\bf 1}[i_0, i_1, \ldots, i_d]|^2$$
where ${\mathbf X}_{\bf 1} = ({\bf 1}\ {\mathbf X})$ is an $m\times (d+1)$ matrix with ${\bf 1}=(1, \ldots, 1)^T$, ${\mathbf X}$ $=$ $({\bf x}_1,$ $\ldots,$ ${\bf x}_m)^T$, and ${\mathbf X}_{\bf 1}[i_0, i_1, \ldots, i_d]$ is the sub-matrix consisting of the $i_0{\rm th}, i_1{\rm th},$ $\ldots,$ $i_d{\rm th}$ rows of ${\mathbf X}_{\bf 1}$~.
\end{lemma}

The proof of Lemma~\ref{coelemma3} is in the supplementary materials (Section~\ref{S3section}). To find D-optimal allocations, we write $|{\mathbf F}|=f(n_1,$ $\ldots,$ $n_m)$ for an order-$(d+J-1)$ homogeneous polynomial function $f$. The {\it D-optimal exact design problem} is to find an integer-valued allocation $(n_1, \ldots, n_m)$ maximizing $f(n_1, \ldots, n_m)$ subject to $n_i \in \{0, 1, \ldots, n\},\ i=1, \ldots, m$ and $n_1 + \cdots + n_m = n$ with given positive integer $n$.
Denote $p_i = n_i/n$, $i=1, \ldots, m$. According to Theorem~\ref{sumFtheorem},
\begin{equation}\label{eq:exactapproximate}
f(n_1, \ldots, n_m)=\left|\sum_{i=1}^m n_i {\mathbf A}_i\right| = \left|n\sum_{i=1}^m p_i {\mathbf A}_i\right| 
= n^{d+J-1} f(p_1, \ldots, p_m)
\end{equation}
Due to \eqref{eq:exactapproximate}, Theorems \ref{polynomialtheorem} and \ref{mindesigntheorem} can be directly applied to approximate design problems too: find a real-valued allocation $(p_1, \ldots, p_m)$ maximizing $f(p_1,p_2,\ldots, p_m)$ subject to $0\leq p_i \leq 1,\ i=1, \ldots, m$ and $p_1 + \cdots + p_m = 1$.

According to Lemma~\ref{ranklemma}, $|{\mathbf A}_{i_03}| > 0$. Thus $c_{\alpha_1, \ldots, \alpha_m}$ in Lemma~\ref{coelemma3} is positive as long as ${\mathbf X}_{\bf 1}[i_0, \ldots, i_d]$ is of full rank. Theorem~\ref{mindesigntheorem} implies that a minimally supported design contains at least $d+1$ support points, while the following theorem states a necessary and sufficient condition for the minimum number of support points to be exactly $d+1$:

\begin{theorem}\label{mindesigntheorem2}
$f({\mathbf p}) > 0$ for some ${\mathbf p}=(p_1, \ldots, p_m)^T$ if and only if the extended design matrix ${\mathbf X}_{\bf 1} = ({\mathbf 1}\ {\mathbf X})$ is of full rank $d+1$.
\end{theorem}

The minimal number of experimental settings required can thus be strictly less than the number of parameters. In the odor removal study, for example, the main-effects cumulative link model (\ref{eq:gamma}) involves four independent parameters -- two $\beta$'s for the covariates ($d=2$) and two $\theta$'s for the intercepts ($J-1=2$) -- while a minimally supported design could involve only three experimental settings. For multinomial responses with $J=3$ categories, we get two degrees of freedom from each experimental setting. Here the optimal allocation of experimental units is often not uniform (see Section~4), contrary to the case of binary responses (\cite{ymm2013, ym2014}).


\fontsize{10.95}{14pt plus.8pt minus .6pt}\selectfont
\setcounter{chapter}{3}
\setcounter{section}{3}
\setcounter{equation}{0} 

\medskip
\noindent {\bf 3. D-optimal Approximate Design}

A (locally) D-optimal approximate design is a real-valued allocation ${\mathbf p}=(p_1, \ldots, p_m)^T$ maximizing $f({\mathbf p}) = f(p_1, \ldots, p_m)$ with pre-specified values of parameters. The solution always exists since $f$ is continuous and the set of feasible allocations
$S:=\{(p_1, \ldots, p_m)^T \in {\mathbb R}^m \mid p_i \geq 0, i=1,\ldots, m;\ \sum_{i=1}^m p_i=1\}$
is convex and compact. A nontrivial D-optimal approximate design problem requires an assumption.

\begin{assumption}\label{assumption3}
$m\geq d+1$ and ${\rm Rank}({\mathbf X}_{\bf 1})=d+1$.
\end{assumption}

Assumption~\ref{assumption3} is adopted throughout. With it, the set of {\it valid} allocations
$S_+:=\{{\bf p}=(p_1, \ldots, p_m)^T\in S\mid f({\bf p}) > 0\}$
is nonempty.
Since ${\mathbf F} =\sum_{i=1}^m n_i {\mathbf A}_i =n\sum_{i=1}^m p_i {\mathbf A}_i$ is linear in ${\mathbf p}$ and $\phi(\cdot)=\log|\cdot|$ is concave on positive semi-definite matrices, $f({\mathbf p}) =n^{1-d-J}|{\mathbf F}| $ is log-concave (\cite{silvey1980}) and thus $S_+$ is also convex.

\begin{theorem}\label{s+theorem}
A feasible allocation ${\mathbf p}=(p_1, \ldots, p_m)^T$ satisfies $f({\mathbf p}) > 0$ if and only if ${\rm Rank}({\mathbf X}_{\bf 1}[\{i\mid p_i>0\}])=d+1$, where ${\mathbf X}_{\bf 1}[\{i\mid p_i>0\}]$ is the sub-matrix consisting of the $\{i\mid p_i>0\}$th rows of ${\mathbf X}_{\bf 1}$~.
\end{theorem}

As a direct conclusion of Theorem~\ref{s+theorem}, $S_+$ contains all ${\mathbf p}$ whose coordinates are all strictly positive. A special case is the uniform allocation ${\mathbf p}_u = (1/m, \ldots, 1/m)^T$.

A necessary and sufficient condition for an approximate design to be D-optimal is of the general-equivalence-theorem type (\cite{kiefer1974, pukelsheim1993, atkinson2007, stufken2012, fedorov2014, ymm2013}), which is convenient when searching for numerical solutions. Following \cite{ymm2013}, for a given ${\mathbf p} = (p_1,$ $\ldots,$ $p_m)^T$ $\in$ $S_+$ and $i \in \{1,\ldots, m\}$, we set
\begin{equation}\label{f_i(z)}
f_i(z) =f\left(\frac{1-z}{1-p_i}p_1,\ldots,\frac{1-z}{1-p_i}p_{i-1},z, \frac{1-z}{1-p_i}p_{i+1},\ldots, \frac{1-z}{1-p_i}p_m\right)
\end{equation}
with $0\leq z\leq 1$. Here $f_i(z)$ is well defined as long as $p_i < 1$.

\begin{theorem}\label{f_i(z)theorem}
Suppose ${\bf p}=(p_1, \ldots, p_m)^T \in S_+$ and $i \in \{1, \ldots, m\}$. For $0\leq z \leq 1$,
\begin{equation}\label{f_i(z)poly}
f_i(z)=(1-z)^d \sum_{j=0}^{J-1} a_j z^j (1-z)^{J-1-j}
\end{equation}
where $a_0=f_i(0)$, $(a_{J-1}, \ldots, a_1)^T={\mathbf B}_{J-1}^{-1}{\mathbf c}$, ${\mathbf B}_{J-1} = (s^{t-1})_{s,t=1, \ldots, J-1}$, and ${\mathbf c}$ $=$ $(c_1, \ldots, c_{J-1})^T$ with $c_j = (j+1)^{d+J-1} j^{-d} f_i(1/(j+1)) - j^{J-1}f_i(0)$.
\end{theorem}

Following the lift-one algorithm proposed in \cite{ymm2013}, we have parallel results and an algorithm for our case. For simplicity, we also call it the {\it lift-one algorithm}.

\begin{theorem}\label{getheorem}
Given an allocation ${\mathbf p} = (p_1^*, \ldots, p_m^*)^T \in S_+$~, ${\mathbf p}$ is D-optimal if and only if for each $i=1, \ldots, m$, $f_i(z), 0\leq z\leq 1$ attains it maximum at $z=p_i^*$~.
\end{theorem}

\noindent A lift-one algorithm
\begin{itemize}
 \item[$1^\circ$] Start with an allocation ${\mathbf p}_0=(p_1,\ldots,p_m)^T$ satisfying $f\left({\mathbf p}_0\right) > 0$.
 \item[$2^\circ$] Set up a random order of $i$ going through $\{1,2,\ldots,m\}$.
 \item[$3^\circ$] For each $i$, determine $f_i(z)$ according to Theorem~\ref{f_i(z)theorem}, with $J$ determinants $f_i(0), f_i(1/2),$ $f_i(1/3),$ $\ldots,$ $f_i(1/J)$ calculated according to (\ref{f_i(z)}).
 \item[$4^\circ$] Use the quasi-Newton method with gradient defined in (\ref{f_i(z)'}) to find $z_*$ maximizing $f_i(z)$ with $0\leq z\leq 1$. If $f_i(z_*) \leq f_i(0)$, let $z_*=0$. Take ${\mathbf p}_*^{(i)} =(p_1(1-z_*)/(1-p_i),$ $\ldots,$ $p_{i-1}(1-z_*)/(1-p_i), z_*, p_{i+1}(1-z_*)/(1-p_i), \ldots, p_m(1-z_*)/(1-p_i))^T$, so $f({\mathbf p}_*^{(i)}) = f_i(z_*)$.
 \item[$5^\circ$] Replace ${\mathbf p}_0$ with ${\mathbf p}_*^{(i)}$, and $f\left({\mathbf p}_0\right)$ with $f({\mathbf p}_*^{(i)})$.
 \item[$6^\circ$] Repeat $2^\circ\sim 5^\circ$ until $f({\mathbf p}_0)=f({\mathbf p}_*^{(i)})$ for each $i$.
\end{itemize}

\begin{theorem}\label{algo1theorem15}
When the lift-one algorithm converges, the resulting ${\mathbf p}$ maximizes $f({\mathbf p})$.
\end{theorem}

\begin{example}\label{odorexample}{\bf Odor removal study}\quad {\rm Here the response was ordinal in nature, {\tt serious odor}, {\tt medium odor}, and {\tt no odor}. We fit the cumulative link model~(\ref{eq:gamma}) to the data presented in Table~\ref{tab:odor}. The estimated values of the model parameters are $(\hat{\beta_1},$ $\hat{\beta_2},$ $\hat{\theta_1},$ $\hat{\theta_2})^T = (-2.44, 1.09, -2.67, -0.21)^T$. If a follow-up experiment is planned and the estimated parameter values are regarded as the true values, the D-optimal approximate allocation found by the lift-one algorithm is ${\mathbf p}_o = (0.4449, 0.2871,$ $0,$ $0.2680)^T$. The efficiency of the uniform ${\mathbf p}_u = (1/4, 1/4,$ $1/4,$ $1/4)^T$ is $(f({\mathbf p}_u)/f({\mathbf p}_o))^{1/4} = 79.7\%$, which is far from satisfactory.
\hfill{$\Box$}}\end{example}

\begin{example}\label{wineexample}{\bf Wine bitterness study}\quad {\rm \citeauthor{christensen2013} (2015, Table~1) aggregated the wine data from \cite{randall1989}. It contains the output of a factorial experiment with two treatment factors each at two levels ({\tt Temperature} $x_1$: cold ($-$) or warm ($+$); {\tt Contact} $x_2$: no ($-$) or yes ($+$)) affecting wine bitterness. The response was ordinal with five levels (from ``1'' being least bitter to ``5'' being most bitter).
The original design employed a uniform allocation ${\mathbf p}_u = (1/4, 1/4, 1/4, 1/4)^T$. The estimated parameter values under the logit link are $(\hat{\beta_1},\hat{\beta_2},\hat{\theta_1},\hat{\theta_2},\hat{\theta_3},\hat{\theta_4})^T = (1.25, 0.76, -3.36, -0.76,  1.45,  2.99)^T$.
If a follow-up experiment is planned regarding the estimated values of the parameters as the true values, then the D-optimal approximate allocation found by the lift-one algorithm is ${\mathbf p}_o = (0.2694, 0.2643, 0.2333, 0.2330)^T$. The efficiency of the original design ${\mathbf p}_u$ is 99.9\%. Nevertheless, the corresponding efficiency may drop to 80\% if $|\beta_1|$ and $|\beta_2|$ are both larger than $3$ (see Figure~1(a)). In that case, the D-optimal allocations are minimally supported, see Figure~1(b); this is discussed further in Section~5.
\hfill{$\Box$}}\end{example}

\begin{figure}[h!]\label{winefig}
\centering
\includegraphics[height=2.5in,width=5in]{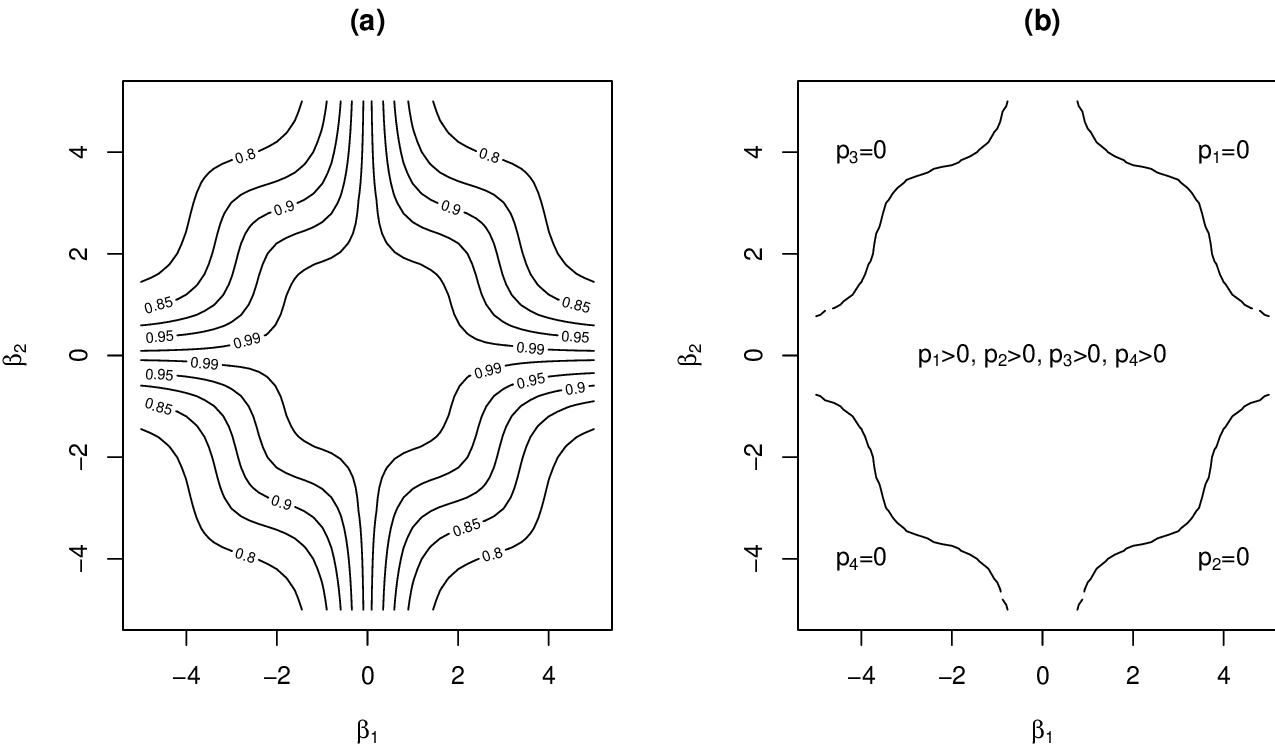}
\footnotesize
\caption{Wine bitterness study with assumed true parameter values $(\beta_1, \beta_2,$ $-3.36,$ $-0.76,$ $1.45, 2.99)^T$: (a) contour plot of efficiency of the original design; (b) regions for a D-optimal design to be minimally supported.}
\normalsize
\end{figure}

In the examples we have studied, the lift-one algorithm often converges within a few iterations. The \cite{ymm2013} lift-one algorithm is guaranteed to converge and can be applied if the lift-one algorithm here does not converge in a pre-specified number of iterations.


\fontsize{10.95}{14pt plus.8pt minus .6pt}\selectfont
\setcounter{chapter}{4}
\setcounter{section}{4}
\setcounter{equation}{0} 

\medskip
\noindent {\bf 4. D-optimal Exact Design}

In the design literature, different discretization methods have been proposed to round an approximate design into an exact design for a given $n$, including the quota method (\cite{kiefer1971, pukelsheim1993}) and the efficient rounding procedure (\cite{pukelsheim1993, pukelsheim1992}), which usually work well for large enough $n$ but with no guarantee for small sample size (\cite{imhofwong2001}).

In this section, we provide a direct search for D-optimal exact designs. From Theorem~\ref{s+theorem}, we have the result as follows:

\begin{corollary}\label{|F|iff}
$|{\mathbf F}| > 0$ if and only if ${\rm Rank}({\mathbf X}_{\bf 1}[\{i\mid n_i>0\}])=d+1$.
\end{corollary}

We assume $n\geq d+1$ throughout this section.
To maximize $f({\mathbf n})=f(n_1,$ $\ldots,$ $n_m)=|{\mathbf F}|$, we adopt the exchange algorithm idea of \cite{fedorov1972}. It is used here to adjust $n_i$ and $n_j$ simultaneously for randomly chosen $(i,j)$, while keeping $n_i+n_j=c$ as a constant.

We start with an ${\mathbf n}=(n_1,\ldots,n_m)^T$ satisfying $f({\mathbf n})>0$.
Following Yang, Mandal, and Majumdar (2016),
for $1\leq i<j\leq m$, let
\begin{equation}\label{f_{ij}}
f_{ij}(z)=
f\left(n_1,\ldots,n_{i-1},z,n_{i+1},\ldots,n_{j-1},c-z,n_{j+1},\ldots,n_m\right)
\end{equation}
where $c=n_i+n_j$, $z=0, 1, \ldots, c$, so $f_{ij}(n_i)=f({\mathbf n})$. From Theorem~\ref{polynomialtheorem}, Lemmas~\ref{coelemma1} and~\ref{coelemma2}, we have the result as follows.

\begin{theorem}\label{f_ij(z)theorem}
Suppose ${\mathbf n}=(n_1,\ldots,n_m)^T$ satisfies $f({\mathbf n})>0$ and $n_i+n_j\geq J$ for given $1\leq i< j \leq m$. For $z=0, 1, \ldots, n_i+n_j$~,
\begin{equation}\label{f_ij(z)poly}
f_{ij}(z)=\sum_{s=0}^J c_s z^s
\end{equation}
where $c_0=f_{ij}(0)$, and $c_1, \ldots, c_J$ can be obtained using $(c_1, \ldots, c_J)^T={\mathbf B}_J^{-1}(d_1,$ $\ldots,$ $d_J)^T$ with ${\mathbf B}_J = (s^{t-1})_{st}$ as a $J \times J$ matrix and $d_s=(f_{ij}(s)-f_{ij}(0))/s$.
\end{theorem}

The $J\times J$ matrix ${\mathbf B}_J$ in Theorem~\ref{f_ij(z)theorem} shares the same form of ${\mathbf B}_{J-1}$ in Theorem~\ref{f_i(z)theorem}. According to Theorem~\ref{f_ij(z)theorem}, in other to maximize $f_{ij}(z)$ with $z=0, 1, \ldots, n_i+n_j$~, one can obtain the exact polynomial form of $f_{ij}(z)$ by calculating $f_{ij}(0), f_{ij}(1), \ldots, f_{ij}(J)$. There is no practical need to find out the exact form of $f_{ij}(z)$ if $n_i + n_j < J$ since one can simply calculate $f_{ij}(z)$ for each $z$. Following \cite{ymm2013}, an exchange algorithm (see the supplementary materials, Section~\ref{S5section}) based on Theorem~\ref{f_ij(z)theorem} could be used to search for a D-optimal exact allocation.

\noindent
{\bf Example~\ref{odorexample} Odor removal study} ({\it continued})\quad
To conduct a follow-up experiment with $n$ experimental units using the exchange algorithm, we obtain the D-optimal exact designs across different $n$'s (Table~\ref{exactdesigns}).
As expected, the D-optimal exact allocation $(n_1, \ldots, n_4)^T$ is consistent with the D-optimal approximate allocation ${\bf p}_o = (p_1, \ldots, p_4)^T$ (last row of Table~\ref{exactdesigns}) for large $n$. The time costs in seconds (last column of Table~\ref{exactdesigns}) are recorded on a PC with 2GHz CPU and 8GB memory. If we rerun an experiment with $n=40$, the D-optimal exact design is ${\mathbf n}_o = (18, 11, 0, 11)^T$, and the efficiency of the uniform design ${\mathbf n}_u = (10, 10, 10, 10)^T$ is $(f({\mathbf n}_u)/f({\mathbf n}_o))^{1/4} = 79.7\%$.
\hfill{$\Box$}

\begin{table}[h!]\caption{D-optimal Exact Designs and the Approximate Design for the Odor Removal Study}\label{exactdesigns}
\begin{center}{\footnotesize
\begin{tabular}{|c|r|r|r|r|r|c|c|}
  \hline
  $n$   & $n_1$ & $n_2$ & $n_3$ & $n_4$ & $n^{-4}|{\mathbf F}|$ & \# iterations & Time(sec.) \\ \hline
  3     & 1 & 1 & 0 & 1 & 0.0002911 & 1 & $<0.01$ \\
  10    & 4 & 3 & 0 & 3 & 0.0003133 & 3 & 0.02 \\
  40    & 18 & 11 & 0 & 11 & 0.0003177 & 3 & 0.02 \\
  100   & 44 & 29 & 0 & 27 & 0.0003180 & 4 & 0.05 \\
  1000  & 445 & 287 & 0 & 268 & 0.0003181 & 5 & 0.39 \\ \hline
  ${\bf p}_o$ & 0.4449 & 0.2871 & 0 & 0.2680 & 0.0003181 & 5 & 0.03 \\
  \hline
\end{tabular}}
\end{center}
\end{table}

\begin{example}\label{industrialexample}{\bf Polysilicon deposition study}\quad {\rm
\cite{wu2008} considered an experiment for studying the polysilicon deposition process with six 3-level factors, described in details by \cite{phadke1989}. Due to the inconvenience of counting the number of surface defects, a major evaluating characteristic, they treated it as a 5-category ordinal variable: 1 for $0\sim 3$ defects, 2 for $4\sim 30$, 3 for $31\sim 300$, 4 for $301\sim 1000$, and 5 for $1001$ and more. The original design, denoted by ${\mathbf n}_u$, includes 18 experimental settings based on an $L_{18}$ orthogonal array. To apply a cumulative link model, we represent each 3-level factor, say $A$, with levels $1,2,3$, by its linear component $A_1$ taking values $-1,0,1$ and a quadratic component $A_2$ taking values $1,-2,1$ (\cite{wu2009}). Then the fitted model with complementary log-log link chosen by both AIC and BIC criteria (see, for example, \cite{agresti2013}) involves four cut-points $(\hat{\alpha_1}, \hat{\alpha_2}, \hat{\alpha_3}, \hat{\alpha_4}) = (-1.59, -0.58,  0.41,  1.22)$, and twelve other coefficients $(\hat{\beta_{11}}, \hat{\beta_{12}}, \hat{\beta_{21}}, \hat{\beta_{22}},$ $\ldots,$ $\hat{\beta_{62}})$ $=$ $(1.45, -0.22,  1.35,  0.02, -0.12, -0.34,  0.19,  0.00,  0.22,$  $0.08,$  $0.05,$  $0.17)$. When the true parameter values were assumed to be the estimated ones, we used the exchange algorithm to find a D-optimal 18-run design, denoted by ${\mathbf n}_o$ (see the supplementary materials, Section~\ref{S6section}, for a list of the 18 experimental settings). Compared with ${\mathbf n}_o$, the efficiency of the original design ${\mathbf n}_u$ is $(f({\mathbf n}_u)/f({\mathbf n}_o))^{1/16} = 73.1\%$.
In order to check the efficiency of a rounded design, we used the lift-one algorithm to find that the D-optimal approximate design contains 100 positive $p_i$'s out of the 729 distinct experimental settings. In this case, both the quota method and the efficient rounding procedure end with the same rounded design ${\mathbf n}_r$ (see Section~\ref{S6section}). Its efficiency is $(f({\mathbf n}_r)/f({\mathbf n}_o))^{1/16} = 86.1\%$.
\hfill{$\Box$}}\end{example}


\fontsize{10.95}{14pt plus.8pt minus .6pt}\selectfont
\setcounter{chapter}{5}
\setcounter{section}{5}
\setcounter{equation}{0} 

\medskip
\noindent {\bf 5. Minimally Supported Design}

It is of practical significance to have experiment run with the minimal number of different settings. For example, the 18 experimental settings in the polysilicon deposition study (Example~\ref{industrialexample}) had to be run in a sequential way and only two settings were arranged on each day (\cite{phadke1989}). Less experimental settings often indicate less time and less cost.
Another practical application of a minimally supported design is that an optimal allocation restricted to those support points can be obtained more easily or even analytically.

According to Theorem~\ref{mindesigntheorem}, a minimally supported design contains at least $d+1$ support points.
On the other hand, according to Theorem~\ref{s+theorem} and Corollary~\ref{|F|iff}, a minimally supported design could contain exactly $d+1$ support points if the extended design matrix ${\mathbf X}_{\bf 1} = ({\bf 1}\ {\mathbf X})$ is of full rank.

\begin{example}\label{d=2example}
{\rm
Let $J=2$ with a binomial response. There are $d+1$ parameters, $\theta_1, \beta_1, \ldots, \beta_d$~. For a general link function $g$ satisfying Assumptions~\ref{pi>0} and \ref{glink}, $g_{i0}=g_{i2}=0$, $g_{i1} = (g^{-1})'(\theta_1-{\mathbf x}_i^T\boldsymbol{\beta})>0$, $e_i = u_{i1} = c_{i1} = g_{i1}^2/[\pi_{i1} (1-\pi_{i1})]$, $i=1, \ldots, m$. Then ${\mathbf A}_{i3}$ in Theorem~\ref{sumFtheorem} contains only the entry $u_{i1}$, and thus $|{\mathbf A}_{i3}|=u_{i1}$, or simply $e_i$ (Lemma~\ref{ranklemma} still holds). Assume further that the $m\times d$ design matrix ${\mathbf X}$ satisfies Assumption~\ref{assumption3}. According to Theorem~\ref{polynomialtheorem}, Lemmas~\ref{coelemma1}, \ref{coelemma2}, and \ref{coelemma3}, given ${\mathbf p}=(p_1, \ldots, p_m)^T$,
\begin{equation}\label{|F|d=2}
f({\mathbf p})=n^{-(d+1)}|{\mathbf F}|=\sum_{1\leq i_0 < i_1 < \cdots < i_d \leq m} |{\mathbf X}_{\bf 1}[i_0, i_1, \ldots, i_d]|^2 p_{i_0} e_{i_0} p_{i_1} e_{i_1} \cdots p_{i_d} e_{i_d}
\end{equation}
Here (\ref{|F|d=2}) is essentially the same as Lemma~3.1 in \cite{ym2014}.
Then a minimally supported design can contain $d+1$ support points and a D-optimal one keeps equal weight $1/(d+1)$ on all support points (\citeauthor{ym2014} (2015, Theorem~3.2)).
\hfill{$\Box$}
}\end{example}

For univariate responses (including binomial ones) under a generalized linear model, a minimally supported design must keep equal weights on all its support points in order to keep D-optimality (\cite{ym2014}). However, for multinomial responses with $J\geq 3$, this is usually not the case. In this section, we use one-predictor ($d=1$) and two-predictor ($d=2$) cases for illustration.

In order to check if a minimally supported design is D-optimal, we need a Karush-Kuhn-Tucker-type condition. Since $f({\mathbf p})$ is log-concave, the Karush-Kuhn-Tucker conditions (\cite{karush1939, kuhn1951}) are also sufficient.

\begin{theorem}\label{KKTtheorem}
An allocation ${\mathbf p} = (p_1^*, \ldots, p_m^*)^T$ satisfying $f({\mathbf p})>0$ is D-optimal if and only if there exists a $\lambda \in \mathbb{R}$ such that $\partial f({\mathbf p})/\partial p_i = \lambda$ if $p_i^* > 0$ or $\leq  \lambda$ if $p_i^* =0$, $i=1, \ldots, m$.
\end{theorem}


\smallskip
\noindent   {\bf 5.1 Minimally supported designs with one predictor}

We start with $d=1$ and $J\geq 3$. The corresponding parameters here are $\beta_1$ and $\theta_1, \ldots, \theta_{J-1}$~.
Consider designs supported on two points ($m=2$, minimally supported), and invoke  Theorem~\ref{polynomialtheorem}, Lemmas~\ref{coelemma1} and \ref{coelemma2}.

\begin{theorem}\label{d=1corollary}
If $d=1$, $J\geq 3$, and $m=2$, the objective function is
\begin{equation}\label{d=1equation}
f(p_1, p_2)=n^{-2}|{\mathbf F}|=\sum_{s=1}^{J-1}c_s p_1^{J-s}p_2^s
\end{equation}
where $(c_1, \ldots, c_{J-1})^T={\mathbf B}_{J-1}^{-1}(d_1, \ldots, d_{J-1})^T$, with ${\mathbf B}_{J-1} = (s^{t-1})_{st}$ as a $(J-1) \times (J-1)$ matrix and $d_s=f(1/(s+1), s/(s+1)) \cdot (s+1)^J/s$.
\end{theorem}

Actually, according to Lemma~\ref{coelemma3}, $c_1 = e_2 \prod_{s=1}^{J-1} g_{1s}^2 \cdot \prod_{t=1}^J \pi_{1t}^{-1} (x_1 - x_2)^2$, $c_{J-1} = e_1 \prod_{s=1}^{J-1} g_{2s}^2 \cdot \prod_{t=1}^J \pi_{2t}^{-1} (x_1 - x_2)^2$, where $x_1, x_2$ are the predictor levels. Theorem~\ref{d=1corollary} provides a way to find the exact form of $f(p_1, p_2)$ after calculating $|{\mathbf F}|$ for $J-1$ different allocations. Then the D-optimal problem is to maximize an order-$J$ polynomial $f(z,1-z)$ for $z \in [0,1]$. As a special case, the D-optimal allocation of $J=3$ can be solved explicitly as follows:

\begin{corollary}\label{d=1J=3corollary}
If $d=1$, $J = 3$, and $m=2$, the objective function is
\begin{equation}\label{d=1J=3equation}
f(p_1, p_2)=p_1 p_2 (c_1 p_1 + c_2 p_2)
\end{equation}
where $c_1 = e_2 g_{11}^2 g_{12}^2 (\pi_{11}\pi_{12}\pi_{13})^{-1} (x_1 - x_2)^2 > 0$, $c_{2} = e_1 g_{21}^2 g_{22}^2 (\pi_{21}\pi_{22}\pi_{23})^{-1}$ $(x_1 - x_2)^2 > 0$, and $x_1, x_2$ are the two levels of the predictor. The D-optimal design ${\mathbf p} = (p_1^*, p_2^*)$ is
\begin{equation}\label{d=1J=3solution}
p_1^* = \frac{c_1-c_2+\sqrt{c_1^2-c_1 c_2 + c_2^2}}{2c_1-c_2+\sqrt{c_1^2-c_1 c_2 + c_2^2}}\ , \quad
p_2^*=\frac{c_1}{2c_1-c_2+\sqrt{c_1^2-c_1 c_2 + c_2^2}}
\end{equation}
Furthermore, $p_1^*=p_2^*=1/2$ if and only if $c_1=c_2$~.
\end{corollary}

Under the setup of Corollary~\ref{d=1J=3corollary}, $p_1^*=p_2^*=1/2$ if $\beta_1=0$. In general $p_1^*\neq p_2^*$, and $p_1^*>p_2^*$ if and only if $c_1 > c_2$.
The following result provides conditions for D-optimality of such a minimally supported design. Its proof is relegated to the supplementary materials (Section~\ref{S3section}).

\begin{corollary}\label{d=1J=3m=3corollary}
Suppose $d=1$, $J = 3$, $m\geq 3$, and let $x_1, \ldots, x_m$ be the $m$ distinct levels of the predictor. A minimally supported design ${\mathbf p} = (p_1^*, p_2^*, 0, \ldots, 0)^T$ is D-optimal if and only if
\begin{itemize}
\item[(1)] $p_1^*, p_2^*$ are defined as in (\ref{d=1J=3solution}),
\item[(2)] $s_{i3} (p_1^*)^2 + (s_{i5} - 2 c_1) p_1^* p_2^* + (s_{i4}-c_2)(p_2^*)^2 \leq 0 $, $i=3, \ldots, m$,
\end{itemize}
where $c_1, c_2$ are as in Corollary~\ref{d=1J=3corollary}, $s_{i3} = e_i g_{11}^2 g_{12}^2 (\pi_{11}\pi_{12}\pi_{13})^{-1} (x_1$ $-$ $x_i)^2 > 0$, $s_{i4} = e_i g_{21}^2 g_{22}^2 (\pi_{21}\pi_{22}\pi_{23})^{-1} (x_2 - x_i)^2 > 0$, $s_{i5} = e_1 (u_{22}u_{i1} + u_{21}u_{i2} - 2b_{22}b_{i2})(x_1 - x_2)(x_1 - x_i) + e_2 (u_{12}u_{i1} + u_{11} u_{i2} - 2 b_{12} b_{i2}) (x_2 - x_1) (x_2 - x_i) + e_i (u_{12} u_{21} + u_{11}u_{22} - 2b_{12} b_{22}) (x_i-x_1)(x_i-x_2)$.
\end{corollary}

\begin{example}\label{d=1J=3example}
{\rm
Consider $d=1$, $J=3$, and $m=3$ with factor levels $\{-1, 0, 1\}$. Under the logit link $g$, the parameters $\beta, \theta_1, \theta_2$ satisfy $g(\gamma_{1j})=\theta_j + \beta$, $g(\gamma_{2j}) = \theta_j$, $g(\gamma_{3j})=\theta_j - \beta$, $j=1,2$. We investigate when a D-optimal design is minimally supported. According to Theorem~\ref{d=1corollary}, a D-optimal deign satisfies $p_1=p_3=1/2$ if $\beta=0$. Figure~2 shows cases with more general parameter values. In Figure~2(a), four regions in $(\theta_1, \theta_2)$-plane are occupied by minimally supported designs ($\theta_1 < \theta_2$ is required). For example, regions labeled with $p_2=0$ indicates a minimally supported design satisfying $p_2=0$ is D-optimal given such a triple $(\theta_1, \theta_2, \beta=-2)$. From Figure~2(b), a design supported on $\{-1, 1\}$ (that is, $p_2=0$) is D-optimal if $\beta$ is not far from $0$.
\hfill{$\Box$}}\end{example}

\begin{figure}[h!]\label{d1J3beta}
\centering
\includegraphics[height=2.5in,width=5in]{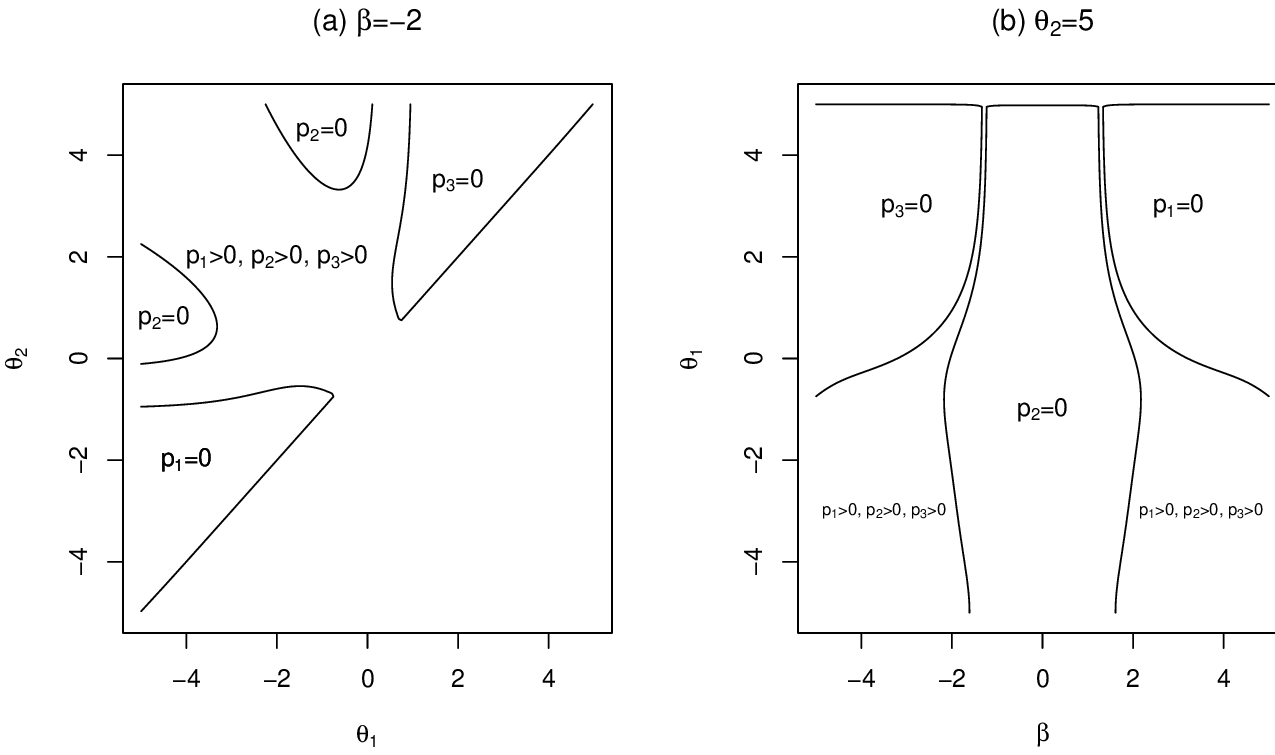}
\caption{Regions for a two-point design to be D-optimal with $d=1$, $J=3$, $x\in\{-1, 0, 1\}$, and logit link (note that $\theta_1 < \theta_2$ is required)}
\end{figure}

\begin{example}\label{d=1J=3m=5example}
{\rm {\bf Toxicity study}\quad
\citeauthor{agresti2013} (2013, Table~8.7) reported data from a developmental toxicity study with one factor (concentration of diEGdiME at five levels: 0, 62.5, 125, 250, 500 mg/kg per day) and a 3-category ordinal response (status of mouse fetus: {\tt nonlive}, {\tt malformation}, or {\tt normal}). In this case, $d=1$, $J=3$, and $m=5$. We fit a cumulative link model with {\tt cauchit} link chosen by both AIC and BIC criteria. The estimated parameter values are $(\hat{\beta_1}, \hat{\theta_1}, \hat{\theta_2})^T = (-0.0176, -8.80, -5.34)$. If $(\hat{\beta_1}, \hat{\theta_1}, \hat{\theta_2})^T$ is regarded as the true parameter value, then the D-optimal approximate allocation found by the lift-one algorithm is ${\mathbf p}_o = (0, 0, 0, 0.4285, 0.5715)^T$, which is minimally supported. Alternatively, for each pair of indices $(i, j)$, $1\leq i < j \leq 5$, we obtain the best design $(p_i^*, p_j^*)$ supported only on $x_i, x_j$ according to Corollary~\ref{d=1J=3corollary}, then check whether $(p_i^*, p_j^*)$ is D-optimal using Corollary~\ref{d=1J=3m=3corollary}. Here ${\mathbf p}_o$ is the only minimally supported design that is also D-optimal. With respect to ${\mathbf p}_o$, the efficiency of the original design (roughly a uniform one) is $52.6\%$.
\hfill{$\Box$}}\end{example}


\noindent   {\bf 5.2 Minimally supported designs with two predictors}

In this section, we consider experiments with two predictors ($d=2$) and a three-category response ($J=3$). The parameters are $\beta_1, \beta_2, \theta_1, \theta_2$~. For cases with $J\geq 4$, similar conclusions could be obtained, but with messier notation.

According to Theorem~\ref{mindesigntheorem}, a minimally supported design needs three support points, for example, $(x_{i1}, x_{i2}),\ i=1,2,3$. Under Assumption~\ref{assumption3}, the $3\times 3$ matrix ${\mathbf X}_{\bf 1} = ({\bf 1}\ {\mathbf X})$ is of full rank. Following Theorem~\ref{polynomialtheorem}, Lemmas~\ref{coelemma1}, \ref{coelemma2}, and \ref{coelemma3}, the objective function with $(d, J, m) = (2, 3, 3)$ is
\begin{equation}\label{d=2m=3equation}
f(p_1, p_2, p_3)=|{\mathbf X}_{\bf 1}|^2 e_1 e_2 e_3 \cdot p_1p_2p_3(w_1p_1 + w_2p_2 + w_3p_3)
\end{equation}
where $w_i=e_i^{-1} g_{i1}^2 g_{i2}^2 (\pi_{i1}\pi_{i2}\pi_{i3})^{-1} > 0$.
Since $f(p_1,$ $p_2,$ $p_3)=0$ if $p_1p_2p_3=0$, we need only consider ${\mathbf p}=(p_1, p_2, p_3)^T$ satisfying $0 < p_1, p_2, p_3 < 1$.

According to Theorem~\ref{KKTtheorem}, ${\mathbf p}$ maximizes $f(p_1, p_2, p_3)$ only if
\begin{equation}\label{3pointcondition}
\frac{\partial f}{\partial p_1} = \frac{\partial f}{\partial p_2} = \frac{\partial f}{\partial p_3}
\end{equation}
Following \cite{tong2014}, we obtain its analytic solution:
\begin{theorem}\label{f4theorem}
Without loss of generality, $w_1 \geq w_2 \geq w_3 > 0$~. The allocation ${\mathbf p}=(p_1^*, p_2^*, p_3^*)^T$ maximizing $f(p_1, p_2, p_3)$ in (\ref{d=2m=3equation}) exists and is unique.  It satisfies $0< p_3^* \leq p_2^* \leq p_1^* < 1$ and can be obtained analytically as follows.
\begin{itemize}
\item[(i)] If $w_1\geq w_2=w_3$, then
$p_1^*=\Delta_1/(4 w_1 + \Delta_1)$, $p_2^*=p_3^*=2 w_1/(4 w_1 + \Delta_1)$,
where $\Delta_1=2 w_1 - 3 w_2 + \sqrt{4 w_1^2 - 4 w_1 w_2 + 9 w_2^2}$~. A special case is $p_1^*$ $= p_2^* = p_3^*=1/3$ if $w_1=w_2=w_3$~.
\item[(ii)] If $w_1=w_2 > w_3$, then
$p_1^*=p_2^*=\Delta_2/[2(\Delta_2+2 w_1)]$, $p_3^*=2w_1/(\Delta_2+2w_1)$,
where $\Delta_2=3 w_1 - 2 w_3 + \sqrt{9 w_1^2 - 4 w_1 w_3 + 4 w_3^2}$~.
\item[(iii)] If $w_1 > w_2 > w_3$, then
$p_1^* = y_1/(y_1 + y_2 + 1)$, $p_2^* = y_2/(y_1 + y_2 + 1)$, $p_3^* = 1/(y_1 + y_2 + 1)$,
where
\[
y_1=-\frac{b_2}{3} - \frac{2^{1/3} (3 b_1 - b_2^2)}{3 A^{1/3}} +  \frac{A^{1/3}}{3\times 2^{1/3}}, \quad
y_2 = \frac{(w_1 - w_3) y_1}{(w_2 - w_3) + (w_1 - w_2) y_1}
\]
with $A=-27 b_0 + 9 b_1 b_2 - 2 b_2^3 + 3^{3/2} (27 b_0^2 + 4 b_1^3 - 18 b_0 b_1 b_2 - b_1^2 b_2^2 + 4 b_0 b_2^3)^{1/2}$, $b_i = c_i/c_3$, $i=0,1,2$, and $c_0 = w_3 (w_2  - w_3) > 0$, $c_1 = 3 w_1 w_2 - w_1 w_3 - 4 w_2 w_3 + 2 w_3^2 > 0$, $c_2 = 2 w_1^2 - 4 w_1 w_2 - w_1 w_3 + 3 w_2 w_3$, $c_3 = w_1 (w_2 -w_1) < 0$.
\end{itemize}
\end{theorem}

The proof of Theorem~\ref{f4theorem} is relegated to the supplementary materials (Section~\ref{S3section}).

\begin{corollary}\label{1/3corollary}
Suppose $d=2$, $J=3$, and $m=3$. Then ${\mathbf p}=(1/3, 1/3,$ $1/3)^T$ is D-optimal if and only if $w_1=w_2=w_3$, where $w_1, w_2, w_3$ are defined as in (\ref{d=2m=3equation}).
\end{corollary}

\begin{example}\label{2^2example}
{\rm
Consider a $2^2$ factorial design problem with a three-category response and four design points $(1, 1)$, $(1, -1)$, $(-1, 1), (-1, -1)$, denoted by $(x_{i1}, x_{i2}), i=1,2,3,4$.
Take $w_i = e_i^{-1} g_{i1}^2 g_{i2}^2 (\pi_{i1} \pi_{i2} \pi_{i3})^{-1},\ i=1,2,3,4$.
There are five special cases: {\it (i)} if $\beta_1=\beta_2=0$, then $w_1=w_2=w_3=w_4$;
{\it (ii)} if $\beta_1=0, \beta_2\neq 0$, then $w_1=w_3$, $w_2=w_4$, but $w_1\neq w_2$;
{\it (iii)} if $\beta_1\neq 0, \beta_2=0$, then $w_1=w_2$, $w_3=w_4$, but $w_1\neq w_3$;
{\it (iv)} if $\beta_1=\beta_2\neq 0$, then $w_2=w_3$, but $w_1, w_2, w_4$ are distinct;
{\it (v)} if $\beta_1=-\beta_2\neq 0$, then $w_1=w_4$, but $w_1, w_2, w_3$ are distinct.
\hfill{$\Box$}}\end{example}

Theorem~\ref{f4theorem} provides analytic forms of minimally supported designs with $d=2$ and $J=3$.

\begin{corollary}\label{d=2J=3m=4corollary}
Suppose $d=2$, $J = 3$, and $m\geq 4$. Let $(x_{i1}, x_{i2}),\ i=1, \ldots, m$ be $m$ distinct level combinations of the two predictors. With ${\mathbf X}_{\bf 1} = ({\bf 1}\ {\mathbf X})$ an $m\times 3$ matrix, a minimally supported design ${\mathbf p} = (p_1^*, p_2^*, p_3^*, 0, \ldots, 0)^T$ is D-optimal if and only if
 $p_1^*, p_2^*, p_3^*$ are obtained according to Theorem~\ref{f4theorem}, and
    \begin{eqnarray*}
    & &|{\mathbf X}_{\bf 1}[1,2,i]|^2 e_1 e_2 e_i p_1^* p_2^* (w_1 p_1^* + w_2 p_2^*) + |{\mathbf X}_{\bf 1}[1,3,i]|^2 e_1 e_3 e_i p_1^* p_3^* (w_1 p_1^* + w_3 p_3^*)\\ &+& |{\mathbf X}_{\bf 1}[2,3,i]|^2 e_2 e_3 e_i p_2^* p_3^* (w_2 p_2^* + w_3 p_3^*) + D_i p_1^* p_2^* p_3^* \\ &\leq & |{\mathbf X}_{\bf 1}[1,2,3]|^2 e_1 e_2 e_3 p_2^* p_3^* (2w_1 p_1^* + w_2 p_2^* + w_3 p_3^*), \hspace{1cm} \mbox{for }i=4, \ldots, m,
    \end{eqnarray*}
    where $e_j = u_{j1} + u_{j2} - 2b_{j2}$, $w_j = e_j^{-1} g_{j1}^2 g_{j2}^2 (\pi_{j1}\pi_{j2}\pi_{j3})^{-1}$, $j=1, \ldots, m$,
    $D_i = \sum_{\{j,k,s,t\} \in E_i} e_j e_k (u_{s1} u_{t2} + u_{s2} u_{t1} - 2 b_{s2} b_{t2}) \cdot |{\mathbf X}_{\bf 1}[j,k,s]|\cdot |{\mathbf X}_{\bf 1}[j,k,t]|$
    with the sum over $E_i = \{(1,2,3,i), (1,3,2,i), (1,i,2,3),$ $(2,3,1,i),$ $(2,i,1,3), (3,i,1,2)\}$.
\end{corollary}

\begin{example}\label{2^2examplelogit}
{\rm
Consider experiments with $d=2$, $J=3$, $m=4$, and design points $(1,\ 1)$, $(1,\ -1)$, $(-1,\ 1)$, $(-1, -1)$.
Figure~3 provides the boundary lines of regions of parameters $(\beta_1, \beta_2, \theta_1, \theta_2)$ for which the best three-point design is D-optimal. In particular, Figure~3(a) shows the region of $(\beta_1, \beta_2)$ for given $\theta_1, \theta_2$~. It clearly indicates that the best three-point design tends to be D-optimal when the absolute values of $\beta_1, \beta_2$ are large. The region tends to be larger as the absolute values of $\theta_1, \theta_2$ increase. On the other hand, Figure~3(b) displays the region of $(\theta_1, \theta_2)$ for given $\beta_1, \beta_2$~. The symmetry of the boundary lines about $\theta_1 + \theta_2=0$ is due to the logit link which is symmetric about $0$. An interesting conclusion based on Corollary~\ref{d=2J=3m=4corollary} is that in this case a three-point design can never be D-optimal if $\beta_1=0$ or $\beta_2=0$.
\hfill{$\Box$}}\end{example}

\begin{figure}\label{beta12}
\centering
\includegraphics[height=2.2in,width=4in]{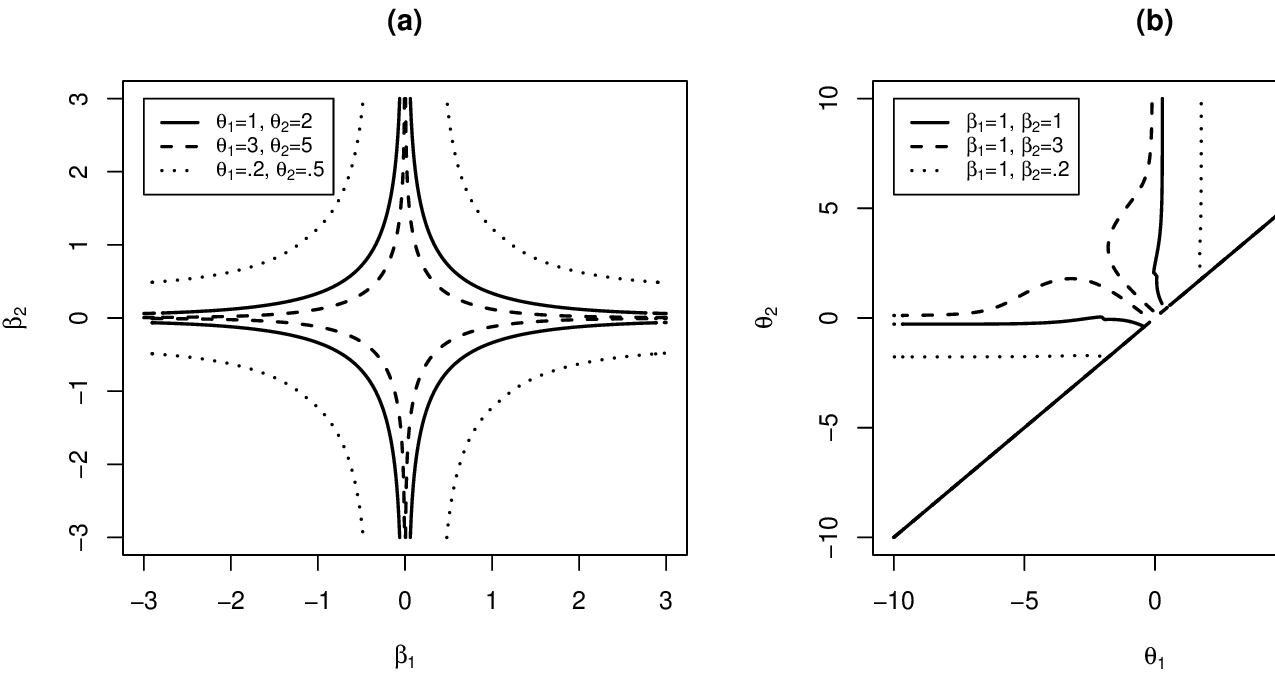}
\caption{Boundary lines for a three-point design to be D-optimal with logit link:
Region of $(\beta_1, \beta_2)$ for given $(\theta_1, \theta_2)$
is outside the boundary lines in Panel (a);
Region of $(\theta_1, \theta_2)$ (with $\theta_1 < \theta_2$) for given
$(\beta_1, \beta_2)$ is between the boundary lines and $\theta_1=\theta_2$ in Panel (b)}
\end{figure}

\begin{remark}\label{rk:dfunif}{\rm
Extra degrees of freedom play an important role against the uniformity of D-optimal allocation in a minimally supported design. For multinomial-type responses with $J$ categories, the total degrees of freedom from $m$ distinct experimental settings is $m(J-1)$, while a cumulative link model contains $d+J-1$ parameters. For a minimally supported design, $m=d+1$ and $m(J-1) = d+J-1$ if and only if $J=2$ (see Example~\ref{d=2example}). Then the objective function $f({\mathbf p}) \propto p_{i_0} p_{i_1} \cdots p_{i_d}$ and the D-optimal allocation is $p_{i_0} = p_{i_1} = \cdots = p_{i_d} = 1/(d+1)$. However, if $J\geq 3$, the degrees of freedom is strictly larger than the number of parameters and there are ``extra" degrees of freedom. In this case, distinct experimental settings may play different roles in estimating the parameters values. For example, if $d=1, J=3, m=2$, the objective function $f({\mathbf p}) = p_1 p_2 (c_1 p_1 + c_2 p_2)$ according to Corollary~\ref{d=1J=3corollary}; if $d=2, J=3, m=3$, $f({\mathbf p}) \propto p_1 p_2 p_3 (w_1 p_1 + w_2 p_2 + w_3 p_3)$ according to equation~\eqref{d=2m=3equation}. The D-optimality of a uniform allocation then depends on $c_1=c_2$ or $w_1=w_2=w_3$, which is not true in general.
}
\end{remark}


\fontsize{10.95}{14pt plus.8pt minus .6pt}\selectfont
\setcounter{chapter}{6}
\setcounter{section}{6}
\setcounter{equation}{0} 

\noindent {\bf 6. EW D-optimal Design}

The previous sections mainly focus on locally D-optimal designs which require assumed parameter values,  $(\beta_1, \ldots,$ $\beta_d,$ $\theta_1,$ $\ldots, \theta_{J-1})$. For many applications, the experimenter may have little information about the values of parameters. Then Bayes D-optimality (\cite{chaloner1995}) which maximizes $E(\log|{\mathbf F}|)$ given a prior distribution on parameters provides a reasonable solution.
An alternative is EW D-optimality (\cite{ymm2013, atkinson2007}) which essentially maximizes $\log|E({\mathbf F})|$. According to \cite{ymm2013}'s simulation study across different models and choices of priors, EW D-optimal designs are much easier to calculate and still highly efficient compared with Bayes designs.

Based on Theorem~\ref{sumFtheorem}, an EW D-optimal design that maximizes $|E({\mathbf F})|$ can be viewed as a locally D-optimal design with $e_i, c_{it}, u_{it}$ and $b_{it}$ replaced by their expectations. After the replacement, Lemma~\ref{cuelemma} still holds. Therefore, almost all results in the previous sections can be applied directly to EW D-optimal designs. The only exception is Lemma~\ref{ranklemma} which provides the formula for $|{\mathbf A}_{i3}|$ in terms of $g_{ij}$ and $\pi_{ij}$~. In order to find EW D-optimal designs, $|{\mathbf A}_{i3}|$ needs to be calculated in terms of $u_{it}$ and $b_{it}$~. For example, $|{\mathbf A}_{i3}|=u_{i1}$ if $J=2$, $|{\mathbf A}_{i3}| = u_{i1} u_{i2} - b_{i2}^2$ if $J=3$, and $|{\mathbf A}_{i3}| = u_{i1} u_{i2} u_{i3} - u_{i1} b_{i3}^2 - u_{i3} b_{i2}^2$ if $J=4$. Then the formulas of $|{\mathbf A}_{i3}|$ in Lemma~\ref{coelemma3}, $c_1, c_2$ in Corollary~\ref{d=1J=3corollary}, $s_{i3}, s_{i4}, s_{i5}$ in Corollary~\ref{d=1J=3m=3corollary}, $w_i$ in (\ref{d=2m=3equation}), and $w_j$ in Corollary~\ref{d=2J=3m=4corollary} need to be written in terms of $u_{it}$ and $b_{it}$.

According to Lemma~\ref{cuelemma}, we only need to calculate $E(u_{it}), i=1, \ldots, m; t=1, \ldots, J-1$ and $E(b_{it}), i=1, \ldots, m; t=2, \ldots, J-1$ (if $J\geq 3$). Then $E(c_{it}) = E(u_{it}) - E(b_{it}) - E(b_{i, t+1})$ and $E(e_i) = \sum_{t=1}^{J-1} E(c_{it})$. After that, we can use the lift-one algorithm in Section~3 or the exchange algorithm in Section~4 to find EW D-optimal designs.

{\bf Example~\ref{odorexample} : Odor Removal Study} ({\it continued})\quad
Instead of assuming the parameter values $(\beta_1, \beta_2, \theta_1,$ $\theta_2)$ $=$ $(-2.44, 1.09, -2.67, -0.21)$, consider true values of parameters that satisfy $\beta_1 \in [-3, -1]$, $\beta_2 \in [0, 2]$, $\theta_1 \in [-4, -2]$, and $\theta_2\in [-1, 1]$. We assume that the four parameters are independently and uniformly distributed within their intervals.
We use {\tt R} function {\tt constrOptim} to maximize $\phi({\mathbf p}) = E(\log |{\mathbf F}|)$ and find the Bayes D-optimal allocation ${\mathbf p}_b = (0.3879, 0.3264, 0.0000, 0.2857)^T$. The procedure costs 313 seconds computational time.
In order to get the EW D-optimal design, we only need 5.43 seconds in total to calculate $E(u_{it})$,  $E(b_{it})$, and find ${\mathbf p}_e = (0.3935, 0.3259,$ $0,$ $0.2806)^T$ using the lift-one algorithm. Even in terms of Bayes Optimality (\cite{chaloner1989,song1998,abebe2014}), the relative efficiency of ${\mathbf p}_e$ with respect to ${\mathbf p}_b$ is $\exp\{(\phi({\mathbf p}_e) - \phi({\mathbf p}_b))/4\} \times 100\% = 99.99\%$, while the relative efficiency of the uniform allocation ${\mathbf p}_u = (0.25, 0.25, 0.25, 0.25)^T$ is $87.67\%$.

In order to check {\it robustness} towards misspecified parameter values, we let $\boldsymbol\theta = (\beta_1, \beta_2,$ $\theta_1,$ $\theta_2)^T$ run through all $0.1$-grid points in $[-3, -1]\times [0, 2]\times [-4, -2] \times [-1, 1]$. For each $\boldsymbol\theta$, we use the lift-one algorithm to find the D-optimal allocation ${\mathbf p}_{\boldsymbol\theta}$ and the corresponding determinant $f({\mathbf p}_{\boldsymbol\theta}) = |{\mathbf F}({\mathbf p}_{\boldsymbol\theta})|$, and then calculate the efficiency $(f({\mathbf p})/f({\mathbf p}_{\boldsymbol\theta}))^{1/4}$ for ${\mathbf p} = {\mathbf p}_b, {\mathbf p}_e$, and ${\mathbf p}_u$, respectively. Table~\ref{R(p,theta)table} shows the summary statistics of the efficiencies. It implies that ${\mathbf p}_b$ and ${\mathbf p}_e$ are comparable and both of them are much better than ${\mathbf p}_u$ in terms of robustness.
\hfill{$\Box$}

\begin{table}[h!]\caption{Summary of Efficiency in Odor Removal Study}\label{R(p,theta)table}
\begin{center}{\footnotesize
\begin{tabular}{|c|r|r|r|r|r|r|}
  \hline
Design & Min. & 1st Quartile & Median  & Mean & 3rd Quartile & Max.\\ \hline
Bayes ${\mathbf p}_b$ & 0.8464 & 0.9813 & 0.9915 & 0.9839 & 0.9964 & 1.0000\\
EW ${\mathbf p}_e$    & 0.8465 & 0.9802 & 0.9917 & 0.9838 & 0.9967 & 1.0000\\
Uniform ${\mathbf p}_u$ & 0.7423 & 0.8105 & 0.8622 & 0.8674 & 0.9249 & 0.9950\\ \hline
\end{tabular}}
\end{center}
\end{table}


\fontsize{10.95}{14pt plus.8pt minus .6pt}\selectfont
\setcounter{chapter}{7}
\setcounter{section}{7}
\setcounter{equation}{0} 

\noindent {\bf 7. Discussion}

In this paper, we use real experiments to illustrate how much improvements the experimenter could make. Compared with our D-optimal designs, the efficiencies of the original designs are often far from satisfactory: 79.7\% in Example~\ref{odorexample}, 73.1\% in Example~\ref{industrialexample}, and 52.6\% in Example~\ref{d=1J=3m=5example}. More interestingly, our D-optimal designs recommended for Example~\ref{odorexample} and Example~\ref{d=1J=3m=5example} are both minimally supported. We have two surprising findings that are different from the cases under univariate generalized linear models \citep{ym2014}: (1) the minimum number of experimental settings can be strictly less than the number of parameters, and (2) the allocation of experimental units on the support points of a minimally supported design is usually not uniform.

Cumulative link models are widely used for modeling ordinal data. Nevertheless, there are other models used for multinomial-type responses, including baseline-category logit model for nominal response, adjacent-categories logit model for ordinal data, and continuation-ratio logit model for hierarchical response (see \cite{agresti2005}, \cite{agresti2013} for a review). The methods developed in this paper could be extended for those models as well. For further extensions, our approaches could be used for planning experiments with more than one categorical response. For example, both the paper feeder experiment and the PCB experiment analyzed by \cite{joseph2004} involved multiple binomial responses.

\vskip 14pt
\noindent {\large\bf Supplementary Materials}

The proofs of Theorems~\ref{sumFtheorem}, \ref{polynomialtheorem}, \ref{mindesigntheorem2}, \ref{s+theorem}, and \ref{f4theorem}, Lemma~\ref{coelemma3}, and Corollaries~\ref{d=1J=3m=3corollary} and \ref{d=2J=3m=4corollary}
are available in the Supplementary Materials.
There are also tabularized formulas for commonly used link functions, additional lemmas for Section~2 and Section~5.2, maximization of $f_i(z)$ in Section 3, exchange algorithm for D-optimal exact allocation in Section 4,
and more results for Example~\ref{industrialexample}.
\par

\vskip 14pt
\noindent {\large\bf Acknowledgements}

We thank Dr.~Suraj Sharma for providing the details of the odor removal study, and Dr.~John Stufken for valuable suggestions on an early version of this paper. We also thank an associate editor and the reviewers for comments and suggestions that substantially improved the quality of the manuscript. This research is in part supported by the LAS Award for Faculty of Science at UIC.
\par

\markboth{\hfill{\footnotesize\rm JIE YANG, LIPING TONG AND ABHYUDAY MANDAL} \hfill}
{\hfill {\footnotesize\rm D-OPTIMAL DESIGNS WITH ORDERED CATEGORICAL DATA} \hfill}

\bibhang=1.7pc
\bibsep=2pt
\fontsize{9}{14pt plus.8pt minus .6pt}\selectfont
\renewcommand\bibname

\vskip .65cm
\noindent
University of Illinois at Chicago
\vskip 2pt
\noindent
E-mail: jyang06@math.uic.edu
\vskip 2pt
\noindent
Phone: (312) 413-3748
\vskip 2pt
\noindent
Fax: (312) 996-1491
\vskip 2pt

\vskip 4pt
\noindent
Advocate Health Care
\vskip 2pt
\noindent
E-mail: lipingtong@gmail.com
\vskip 2pt

\vskip 4pt
\noindent
University of Georgia
\vskip 2pt
\noindent
E-mail: amandal@stat.uga.edu




\clearpage
\setcounter{page}{1}
\def\thepage{S\arabic{page}}


\markboth{\hfill{\footnotesize\rm JIE YANG, LIPING TONG AND ABHYUDAY MANDAL} \hfill}
{\hfill {\footnotesize\rm SUPPLEMENTARY MATERIALS} \hfill}

\fontsize{10.95}{14pt plus.8pt minus .6pt}\selectfont
\vspace{0.8pc}
\centerline{\large\bf D-OPTIMAL DESIGNS WITH}
\vspace{2pt} \centerline{\large\bf ORDERED CATEGORICAL DATA}
\vspace{.25cm}
\centerline{Jie Yang$^{1}$, Liping Tong$^{2}$ and Abhyuday Mandal$^{3}$ }
\vspace{.4cm}
\centerline{\it  $^1$University of Illinois at Chicago, $^2$Advocate Health Care and $^3$University of Georgia}
\vspace{.55cm}
 \centerline{\bf Supplementary Materials}
\vspace{.55cm}
\fontsize{9}{11.5pt plus.8pt minus .6pt}\selectfont
\par

\setcounter{section}{0}
\setcounter{equation}{0}
\setcounter{table}{0}
\setcounter{lemma}{0}
\def\theequation{S.\arabic{equation}}
\def\thelemma{S.\arabic{lemma}}
\def\thetable{S.\arabic{table}}
\def\thesection{S}

\fontsize{12}{14pt plus.8pt minus .6pt}\selectfont

\subsection{Commonly Used Link Functions for Cumulative Link Models}\label{S1section}
\begin{center}
\begin{tabular}{|l|c|c|c|}\hline
Link function         & \quad $g(\gamma)$\quad\quad                & \quad $g^{-1}(\eta)$\quad\quad               & \quad $(g^{-1})'(\eta)$\quad\quad  \\
\hline
logit                 & $\log\left(\frac{\gamma}{1-\gamma}\right)$ & $\frac{e^\eta}{1+e^\eta}$                    & $\frac{e^\eta}{(1+e^\eta)^2}$      \\
probit                & $\Phi^{-1}(\gamma)$                        & $\Phi(\eta)$                                 & $\phi(\eta)$                       \\
log-log               & $-\log[-\log (\gamma)]$                    & $\exp\{-e^{-\eta}\}$                         & $\exp\{-\eta-e^{-\eta}\}$          \\
c-log-log & $\log[-\log (1-\gamma)]$                   & $1-\exp\{-e^\eta\}$                          & $\exp\{\eta - e^\eta\}$            \\
cauchit               & $\tan[\pi(\gamma-\frac{1}{2})]$            & $\frac{1}{\pi} \arctan(\eta) + \frac{1}{2}$  & $\frac{1}{\pi(1+\eta^2)}$          \\
\hline
\end{tabular}
\end{center}
where $\Phi^{-1}(\cdot)$ is the cumulative distribution function of $N(0,1)$, $\phi(\cdot)$ is the probability density function of $N(0,1)$, and ``c-log-log" stands for complementary log-log.

\noindent {\bf Example~\ref{logitd2J3example}} ({\it continued}) \quad For logit link $g$, $g^{-1}(\eta) = e^\eta/(1+e^\eta)$ and $(g^{-1})' = g^{-1}(1-g^{-1})$. Thus
$g_{ij} = (g^{-1})'(\theta_j-{\mathbf x}_i^T\boldsymbol{\beta})=\gamma_{ij}(1-\gamma_{ij})$. With $J=3$, we have $\pi_{i1} + \pi_{i2} + \pi_{i3} = 1$ for $i=1, \ldots, m$. Then for $i=1, \ldots, m$,
$
g_{i1} = \pi_{i1} (\pi_{i2} + \pi_{i3}),
g_{i2} = (\pi_{i1} + \pi_{i2}) \pi_{i3},
b_{i2} = \pi_{i1}\pi_{i3}\pi_{i2}^{-1} (\pi_{i1} + \pi_{i2}) (\pi_{i2} + \pi_{i3}),
u_{i1} = \pi_{i1}\pi_{i2}^{-1}(\pi_{i1}+\pi_{i2})(\pi_{i2}+\pi_{i3})^2,
u_{i2} = \pi_{i3}\pi_{i2}^{-1}(\pi_{i1}+\pi_{i2})^2(\pi_{i2}+\pi_{i3}),
c_{i1} = \pi_{i1} (\pi_{i1}+\pi_{i2}) (\pi_{i2}+\pi_{i3}),
c_{i2} = \pi_{i3} (\pi_{i1}+\pi_{i2}) (\pi_{i2}+\pi_{i3}),
e_i = (\pi_{i1}+\pi_{i2})(\pi_{i1}+\pi_{i3})(\pi_{i2}+\pi_{i3})
$.
\hfill{$\Box$}

\subsection{Additional Lemmas}\label{S2section}

{\bf For Section~2:} Since $(Y_{i1}, \ldots, Y_{iJ}), i=1, \ldots, m$ are $m$ independent random vectors, the log-likelihood function (up to a constant) of the cumulative link model is
$$l(\beta_1, \ldots, \beta_d, \theta_1, \ldots,\theta_{J-1}) = \sum_{i=1}^m \sum_{j=1}^J Y_{ij}\log (\pi_{ij})$$
where $\pi_{ij} = \gamma_{ij} - \gamma_{i,j-1}$ with $\gamma_{ij} = g^{-1}(\theta_j - {\mathbf x}_i^T \boldsymbol{\beta})$ for $j=1,\ldots, J-1$ and $\gamma_{i0}=0$, $\gamma_{iJ}=1$, $i=1, \ldots, m$.
For $s=1, \ldots, d$, $\ t=1, \ldots, J-1$,
\begin{eqnarray*}
\frac{\partial l}{\partial \beta_s} &=& \sum_{i=1}^m (-x_{is})\cdot \left\{
\frac{Y_{i1}}{\pi_{i1}}\cdot (g^{-1})'(\theta_1 - {\mathbf x}_i^T\boldsymbol{\beta})\right.\\
& & + \frac{Y_{i2}}{\pi_{i2}}\cdot\left[(g^{-1})'(\theta_2 - {\mathbf x}_i^T\boldsymbol{\beta})
- (g^{-1})'(\theta_1 - {\mathbf x}_i^T\boldsymbol{\beta})\right]\\
& & + \cdots
\left.+ \frac{Y_{iJ}}{\pi_{iJ}}\left[-(g^{-1})'(\theta_{J-1} - {\mathbf x}_i^T\boldsymbol{\beta})\right]\right\}\\
\frac{\partial l}{\partial \theta_t} &=& \sum_{i=1}^m (g^{-1})'(\theta_t - {\mathbf x}_i^T\boldsymbol{\beta}) \left(\frac{Y_{it}}{\pi_{it}}-\frac{Y_{i,t+1}}{\pi_{i,t+1}}\right)
\end{eqnarray*}

Since $Y_{ij}$'s come from multinomial distributions, we know $E(Y_{ij})=n_i\pi_{ij}$~, $E(Y_{ij}^2) = n_i(n_i-1)\pi_{ij}^2 + n_i \pi_{ij}$~, and $E(Y_{is}Y_{it})=n_i(n_i-1)\pi_{is}\pi_{it}$ when $s\neq t$.
Then we have the following lemma:

\begin{lemma}\label{Fstlemma}
Let ${\mathbf F}=(F_{st})$ be the $(d+J-1)\times (d+J-1)$ Fisher information matrix.
\begin{itemize}
\item[(i)] For $1\leq s\leq d$, $1\leq t\leq d$,
$$F_{st}=E\left(\frac{\partial l}{\partial \beta_s}\frac{\partial l}{\partial \beta_t}\right)
=\sum_{i=1}^m n_i x_{is}x_{it}\sum_{j=1}^J \frac{(g_{ij}-g_{i,j-1})^2}{\pi_{ij}}$$
where $g_{ij}=(g^{-1})'(\theta_j-{\mathbf x}_i^T\boldsymbol{\beta})>0$ for $j=1, \ldots, J-1$ and $g_{i0}=g_{iJ}=0$.
\item[(ii)] For $1\leq s\leq d$, $1\leq t\leq J-1$,
$$F_{s,d+t}=E\left(\frac{\partial l}{\partial \beta_s}\frac{\partial l}{\partial \theta_t}\right)
=\sum_{i=1}^m n_i (-x_{is})g_{it}\left(\frac{g_{it}-g_{i,t-1}}{\pi_{it}} - \frac{g_{i,t+1}-g_{it}}{\pi_{i,t+1}}\right)$$
\item[(iii)] For $1\leq s\leq J-1$, $1\leq t\leq d$,
$$F_{d+s,t}=E\left(\frac{\partial l}{\partial \theta_s}\frac{\partial l}{\partial \beta_t}\right)
=\sum_{i=1}^m n_i (-x_{it})g_{is}\left(\frac{g_{is}-g_{i,s-1}}{\pi_{is}} - \frac{g_{i,s+1}-g_{is}}{\pi_{i,s+1}}\right)$$
\item[(iv)] For $1\leq s\leq J-1$, $1\leq t\leq J-1$,
\[
F_{d+s, d+t}=E\left(\frac{\partial l}{\partial \theta_s}\frac{\partial l}{\partial \theta_t}\right)
=\left\{\begin{array}{ll}
\sum_{i=1}^m n_i g_{is}^2(\pi_{is}^{-1}+\pi_{i,s+1}^{-1}), & \mbox{ if }s=t\\
\sum_{i=1}^m n_i g_{is}g_{it}(-\pi_{i,s\vee t}^{-1}), & \mbox{ if }|s-t|=1\\
0, & \mbox{ if }|s-t|\geq 2
\end{array}
\right.
\]
where $s\vee t=\max\{s,t\}$.
\end{itemize}
\end{lemma}

\cite{perevozskaya2003} obtained a detailed form of Fisher information matrix for logit link and one predictor. Our expressions here are good for fairly general link and $d$ predictors.
To simplify the notations, we denote for $i=1, \ldots, m$,
\begin{eqnarray}
e_i &=& \sum_{j=1}^J \frac{(g_{ij}-g_{i,j-1})^2}{\pi_{ij}} > 0\label{e_i}\\
c_{it} &=& g_{it}\left(\frac{g_{it}-g_{i,t-1}}{\pi_{it}} - \frac{g_{i,t+1}-g_{it}}{\pi_{i,t+1}}\right), \quad t=1, \ldots, J-1\label{c_it}\\
u_{it} &=& g_{it}^2(\pi_{it}^{-1}+\pi_{i,t+1}^{-1}) > 0, \quad t=1, \ldots, J-1\label{u_it}\\
b_{it} &=& g_{i,t-1}g_{it}\pi_{it}^{-1} > 0, \quad t=2, \ldots, J-1\ (\mbox{if }J\geq 3)\label{b_ij}
\end{eqnarray}
Note that $g_{ij}$ is defined in Lemma~\ref{Fstlemma} (i). Then we obtain the following lemma which plays a key role in calculating $|{\mathbf F}|$.

\begin{lemma}\label{cuelemma}
$c_{it} = u_{it} - b_{it} - b_{i, t+1}$, $i=1, \ldots, m; t=1, \ldots, J-1$;
$e_i = \sum_{t=1}^{J-1} c_{it}=\sum_{t=1}^{J-1} (u_{it} - 2 b_{it})$,  $i=1, \ldots, m$,
where $b_{i1}=b_{iJ}=0$ for $i=1, \ldots, m$.
\end{lemma}

\begin{lemma}\label{rank1lemma}
${\rm Rank}(({\mathbf A}_{i1}\> {\mathbf A}_{i2})) \leq 1$ where ``=" is true if and only if ${\mathbf x}_i\neq 0$.
\end{lemma}

Based on Lemmas~\ref{ranklemma} and \ref{rank1lemma}, we obtain the two lemmas below on $c_{\alpha_1, \ldots, \alpha_m}$ which significantly simplify the structure of $|{\mathbf F}|$ as a polynomial of $(n_1, \ldots, n_m)$.

\begin{lemma}\label{coelemma1}
If $\max_{1\leq i\leq m}\alpha_i \geq J$, then $|{\mathbf A}_\tau|=0$ for any $\tau\in (\alpha_1, \ldots, \alpha_m)$ and thus $c_{\alpha_1, \ldots, \alpha_m} =0$.
\end{lemma}

\noindent {\bf Proof of Lemma~\ref{coelemma1}:} Without any loss of generality, we assume
$\alpha_1\geq \alpha_2\geq \cdots \geq \alpha_m$~. Then $\max_{1\leq i\leq m}\alpha_i \geq J$ implies $\alpha_1 \geq J$. In this case, for any $\tau \in (\alpha_1, \ldots, \alpha_m)$, $\tau^{-1}(1) := \{i \mid \tau(i)=1\} $ $\subset$ $\{1, \ldots, d+J-1\}$ and $|\tau^{-1}(1)|=\alpha_1$~. If $|\tau^{-1}(1) \cap \{1, \ldots, d\}| \geq 2$, then $|A_\tau|=0$ due to Lemma~\ref{rank1lemma}; otherwise $\{d+1, \ldots, d+J-1\} \subset \tau^{-1}(1)$ and thus $|A_\tau|=0$ due to Lemma~\ref{ranklemma}. Thus $c_{\alpha_1, \ldots, \alpha_m} = 0$ according to (\ref{Atau}) provided in Theorem~\ref{polynomialtheorem}.
\hfill{$\Box$}

\begin{lemma}\label{coelemma2}
If $\#\{i:\alpha_i\geq 1\} \leq d$, then $|{\mathbf A}_\tau|=0$ for any $\tau\in (\alpha_1, \ldots, \alpha_m)$ and thus $c_{\alpha_1, \ldots, \alpha_m} =0$.
\end{lemma}

\noindent {\bf Proof of Lemma~\ref{coelemma2}:} Without any loss of generality, we assume
$\alpha_1\geq \alpha_2\geq \cdots \geq \alpha_m$~. Then $\#\{i:\alpha_i\geq 1\} \leq d$ indicates $\alpha_{d+1} = \cdots = \alpha_m = 0$.
Let $\tau : \{1, 2, \ldots, d+J-1\} \rightarrow \{1, \ldots, m\}$ satisfy $\tau \in (\alpha_1, \ldots, \alpha_m)$. Then the $(d+J-1)\times (d+J-1)$ matrix $A_\tau$ can be written as
\[
\left(
\begin{array}{cc}
A_{\tau 1} & A_{\tau 2} \\
A_{\tau 3} & A_{\tau 4}\\
\end{array}
\right)
\]
\[
=
\left(
\begin{array}{cc}
(e_{\tau(s)} x_{\tau(s)s}x_{\tau(s)t})_{s=1,\ldots d; t=1,\ldots,d} & (-x_{\tau(s)s}c_{\tau(s)t})_{s=1,\ldots,d; t=1,\ldots,J-1} \\
(-c_{\tau(d+s)s}x_{\tau(d+s)t})_{s=1,\ldots,J-1; t=1,\ldots,d} & A_{\tau 4}\\
\end{array}
\right)
\]
where the $(J-1)\times (J-1)$ matrix $A_{\tau 4}$ is either a single entry $u_{\tau(d+1)1}$ (if $J=2$) or symmetric tri-diagonal with diagonal entries $u_{\tau(d+1)1}, \ldots, u_{\tau(d+J-1), J-1}$,  upper off-diagonal entries $-b_{\tau(d+1)2}, \ldots, -b_{\tau(d+J-2),J-1}$, and lower off-diagonal entries $-b_{\tau(d+2)2},$ $\ldots,$ $-b_{\tau(d+J-1),J-1}$~.
Note that $A_\tau$ is asymmetric in general.

If $\#\{i:\alpha_i\geq 1\} \leq d-1$, then there exists an $i_0$ such that $1\leq i_0 \leq d$ and $|\tau^{-1}(i_0) \cap \{1, \ldots, d\}|\geq 2$. In this case, $|A_\tau|=0$ according to Lemma~\ref{rank1lemma}.

If $\#\{i:\alpha_i\geq 1\} = d$, we may assume $|\tau^{-1}(i) \cap \{1, \ldots, d\}| = 1$ for $i=1, \ldots, d$ (otherwise $|A_\tau|=0$ according to Lemma~\ref{rank1lemma}).
Suppose $\alpha_1 \geq \alpha_2 \geq \cdots \geq \alpha_k \geq 2 > \alpha_{k+1}$~. Then $\{d+1, \ldots, d+J-1\} \subset \cup_{i=1}^k \tau^{-1}(i)$ and $\sum_{i=1}^k (\alpha_i-1) = J-1$.
In order to show $|A_\tau|=0$, we first replace $A_{\tau 1}$ with $A^{(1)}_{\tau 1} = (e_{\tau(s)} x_{\tau(s)t})_{s=1,\ldots d;\ t=1,\ldots,d}$ and replace $A_{\tau 2}$ with $A^{(1)}_{\tau 2} = (-c_{\tau(s)t})_{s=1,\ldots,d;\ t=1,\ldots,J-1}$. It changes $A_\tau$ into a new matrix $A^{(1)}_\tau$~.
Note that $|A_\tau|=\prod_{s=1}^dx_{\tau(s)s}\cdot |A^{(1)}_\tau|$. According to Lemma~\ref{cuelemma}, the sum of the columns of $A^{(1)}_{\tau 2}$ is $(-e_{\tau(1)}, \ldots, -e_{\tau(d)})^T$, and the elementwise sum of the columns of $A_{\tau 4}$ is $(c_{\tau(d+1)1}, c_{\tau(d+2)2},$ $\ldots,$ $c_{\tau(d+J-1),J-1})^T$.
Secondly, for $t=1, \ldots, d$, we add $x_{1t}(-e_{\tau(1)}, \ldots, -e_{\tau(d)}, c_{\tau(d+1)1}, $ $\ldots,$ $c_{\tau(d+J-1),J-1})^T$ to the $t$th column of $A^{(1)}_\tau$. We denote the resulting matrix by $A^{(2)}_\tau$. Note that $|A^{(1)}_\tau| = |A^{(2)}_\tau|$. We consider the sub-matrix $A^{(2)}_{\tau d}$ which consists of the first $d$ columns of $A^{(2)}_\tau$. For $s \in \tau^{-1}(1)$, the $s$th row of $A^{(2)}_{\tau d}$ is simply $0$. For $i=2, \ldots, k$, the $j$th row of $A^{(2)}_{\tau d}$ is proportional to $(x_{i1}-x_{11}, x_{i2} - x_{12}, \ldots, x_{id} - x_{1d})$ if $j \in \tau^{-1}(i)$. Therefore, ${\rm Rank}(A^{(2)}_{\tau d}) \leq (d+J-1) - \alpha_1 - \sum_{i=2}^k (\alpha_i - 1) = d - 1$, which leads to $|A^{(2)}_\tau|=0$ and thus $|A^{(1)}_\tau|=0$, $|A_\tau|=0$. According to (\ref{Atau}) in Theorem~\ref{polynomialtheorem}, $c_{\alpha_1, \ldots, \alpha_m} = 0$.
\hfill{$\Box$}

\begin{lemma}\label{|F|lemma}
${\mathbf F}={\mathbf F}({\mathbf p})$ is always positive semi-definite. It is positive definite if and only if ${\mathbf p}\in S_+$~. Furthermore, $\log f({\mathbf p})$ is concave on $S$.
\end{lemma}

\bigskip
\noindent
{\bf For Section 5.2:} The procedure seeking for analytic solutions here follows \cite{tong2014}.
As a direct conclusion of the Karush-Kuhn-Tucker conditions (see also Theorem~\ref{KKTtheorem}), a necessary condition for $(p_1, p_2, p_3)$ to maximize $f(p_1, p_2, p_3)$ in (\ref{d=2m=3equation}) is (\ref{3pointcondition}), which are equivalent to $\partial f/\partial p_1 = \partial f/\partial p_3$ and $\partial f/\partial p_2 = \partial f/\partial p_3$~. In terms of $p_i, w_i$'s, they are
\begin{eqnarray}
(p_3 - p_1) (p_1 w_1 + p_2 w_2 + p_3 w_3) &=&  (w_3 - w_1) p_1 p_3\label{f41eq} \\
(p_3 - p_2) (p_1 w_1 + p_2 w_2 + p_3 w_3) &=&  (w_3 - w_2) p_2 p_3\label{f42eq}
\end{eqnarray}
Denote $y_1=p_1/p_3 > 0$ and $y_2=p_2/p_3>0$. Since $p_1 + p_2 + p_3 = 1$, it implies $p_3 = 1/(y_1 + y_2 + 1)$, $p_1 = y_1/(y_1 + y_2 + 1)$, and $p_2 = y_2/(y_1 + y_2 + 1)$. In terms of $y_1, y_2$, (\ref{f41eq}) and (\ref{f42eq}) are equivalent to
\begin{eqnarray}
(1 - y_1) (y_1 w_1 + y_2 w_2 + w_3) &=&  (w_3 - w_1) y_1\label{f43eq} \\
(1 - y_2) (y_1 w_1 + y_2 w_2 + w_3) &=&  (w_3 - w_2) y_2\label{f44eq}
\end{eqnarray}

\begin{lemma}\label{monolemma}
Suppose $0<w_3 < w_2 < w_1$~. If $(p_1, p_2, p_3)$ maximizes $f(p_1, p_2, p_3)$ in (\ref{d=2m=3equation}) under the constrains $p_1, p_2, p_3\geq 0$ and $p_1+p_2+p_3=1$, then $0<p_3\leq p_2\leq p_1<1$.
\end{lemma}

The proof of the lemma above is straightforward, because otherwise one could exchange $p_i, p_j$ to strictly improve $f(p_1, p_2, p_3)$~. Now we are ready to get solutions to equations (\ref{f43eq}) and (\ref{f44eq}) case by case.
\begin{itemize}
\item[(i)] $w_1=w_3$~. In that case, (\ref{f43eq}) implies $y_1=1$. After plugging it into (\ref{f44eq}), the only positive solution is\\
$y_2=(-3 w_1 + 2 w_2 + \sqrt{9 w_1^2 - 4 w_1 w_2 + 4 w_2^2})/(2 w_2)$~.

\item[(ii)] $w_2=w_3$~. In that case, (\ref{f44eq}) implies $y_2=1$. After plugging it into (\ref{f43eq}), the only positive solution is\\
$y_1=(2 w_1 - 3 w_2 + \sqrt{4 w_1^2 - 4 w_1 w_2 + 9 w_2^2})/(2 w_1)$~.

\item[(iii)] $w_1=w_2$ but $w_1\neq w_3$~. The ratio of (\ref{f43eq}) and (\ref{f44eq}) leads to $y_1=y_2$~. After plugging it into (\ref{f43eq}), the only positive solution is
$y_1=(3 w_1 - 2 w_3 + \sqrt{9 w_1^2 - 4 w_1 w_3 + 4 w_3^2})/(4 w_1)$~.

\item[(iv)] $w_1, w_2, w_3$ are distinct. Without any loss of generality, we assume $0<w_3 < w_2 < w_1$,  because otherwise the previous elimination procedure in the order of $p_3, p_2, p_1$ could be easily changed accordingly. Based on Lemma~\ref{monolemma}, if $(p_1, p_2, p_3)$ maximizes $f_4$, then $0<p_3\leq p_2\leq p_1<1$ and thus $y_1\geq y_2 \geq 1$.
The ratio of (\ref{f43eq}) and (\ref{f44eq}) leads to $(1-y_1)/(1-y_2)=(w_3-w_1)/(w_3-w_2)\cdot y_1/y_2$, which implies
\begin{equation}\label{f4y2}
y_2 = \frac{(w_1 - w_3) y_1}{(w_2 - w_3) + (w_1 - w_2) y_1}~.
\end{equation}
Note that $(w_2 - w_3) + (w_1 - w_2) y_1\geq w_1 - w_3 > 0$. After plugging (\ref{f4y2}) into (\ref{f43eq}), we get
\begin{equation}\label{y1equation}
c_0 + c_1 y_1 + c_2 y_1^2 + c_3y_1^3=0
\end{equation}
where
$c_0 = w_3 (w_2  - w_3) > 0$,
$c_1 = 3 w_1 w_2 - w_1 w_3 - 4 w_2 w_3 + 2 w_3^2 > 0$,
$c_2 = 2 w_1^2 - 4 w_1 w_2 - w_1 w_3 + 3 w_2 w_3$,
$c_3 = w_1 (w_2 -w_1) < 0$.
\end{itemize}

\begin{lemma}\label{y1lemma}
Suppose $0<w_3 < w_2 < w_1$~. Then equation~(\ref{y1equation}) has one and only one solution $y_1^* \geq 1$. Furthermore, $y_1^* > 1$.
\end{lemma}

\noindent {\bf Proof of Lemma~\ref{y1lemma}:}\quad
In order to locate the roots of equation~(\ref{y1equation}), we let $f_1(y_1)=c_0 + c_1 y_1 + c_2 y_1^2 + c_3y_1^3$. Then $f_1(1)=c_0+c_1+c_2+c_3$ $=$ $(w_1 - w_3)^2 > 0$.

On the other hand, the first derivative of $f_1$ is
$f'_1(y_1) = a_0 + a_1 y_1 + a_2 y_1^2$~,
where $a_0=3 w_1 w_2 - w_1 w_3 - 4 w_2 w_3 + 2 w_3^2 = w_1 (w_2 - w_3) + 2 (w_1 - w_2) w_2 + 2 (w_2 - w_3)^2 > 0$, $a_1 =  2 (2 w_1^2 - 4 w_1 w_2 - w_1 w_3 + 3 w_2 w_3)$, and $a_2 = 3 w_1 (w_2 - w_1) < 0$.
Therefore, $a_1^2 - 4 a_0 a_2 > a_1^2 \geq 0$ and
$f'_1(y_1)=a_2 (y_1 - y_{11}) (y_1 - y_{12})$,
where
$$y_{11} = \frac{-a_1 + \sqrt{a_1^2 - 4 a_0 a_2}}{2 a_2} < 0,\quad
y_{12} = \frac{-a_1 - \sqrt{a_1^2 - 4 a_0 a_2}}{2 a_2} > y_{11}
$$
It can be verified that $y_{12} < 1$ if and only if $w_1 < 2(w_2+w_3)$. There are two cases: {\it Case  (i)}: If $y_{12} < 1$, then $f_1'(y_1) < 0$ for all $y_1 > 1$. That is, $f_1(y_1)$ strictly decreases after $y_1=1$. Since  $f_1(1) > 0$ and $f_1(\infty)=-\infty$, then there is one and only one solution in $(1, \infty)$; {\it Case (ii)}: If $y_{12} \geq 1$, then $f_1'(y_1)\geq 0$ for $y_1 \in [1, y_{12}]$ and $f_1'(y_1) < 0$ for $y_1 \in (y_{12}, \infty)$. That is, $f_1(y_1)$ increases in $[1, y_{12}]$ and then strictly decreases in $(y_{12}, \infty)$. Again, due to $f_1(1) > 0$ and $f_1(\infty)=-\infty$, there is one and only one solution in $(1, \infty)$. In either case, the conclusion is justified.
\hfill{$\Box$}

\subsection{Additional Proofs}\label{S3section}

\noindent {\bf Proof of Theorem~\ref{sumFtheorem}}\quad
It is a direct conclusion of Lemmas~\ref{Fstlemma} and \ref{cuelemma}.
\hfill{$\Box$}

Examples of ${\mathbf A}_{i3}$ in Theorem~\ref{sumFtheorem} include $(u_{i1})$,
\[
\left(
\begin{array}{rr}
u_{i1} & -b_{i2}\\
-b_{i2} & u_{i2}
\end{array}
\right),
\left(
\begin{array}{rrr}
u_{i1}  & -b_{i2} & 0\\
-b_{i2} & u_{i2}  & -b_{i3}\\
0       & -b_{i3} & u_{i3}
\end{array}
\right),
\left(
\begin{array}{rrrr}
u_{i1}  & -b_{i2} & 0       & 0\\
-b_{i2} & u_{i2}  & -b_{i3} & 0\\
0       & -b_{i3} & u_{i3}  & -b_{i4}\\
0       & 0       & -b_{i4} & u_{i4}
\end{array}
\right)
\]
for $J=2$, $3$, $4$, or $5$ respectively.

\bigskip\noindent
{\bf Proof of Theorem~\ref{polynomialtheorem}}\quad
To study the structure of $|{\mathbf F}|$ as a polynomial function of $(n_1, \ldots, n_m)$, we denote the $(k,l)$th entry of ${\mathbf A}_i$ by $a_{kl}^{(i)}$. Given a row map $\tau : \{1, 2, \ldots, d+J-1\} \rightarrow \{1, \ldots, m\}$, we define a $(d+J-1)\times (d+J-1)$ matrix ${\mathbf A}_\tau = \left(a_{kl}^{(\tau(k))}\right)$ whose $k$th row is given by the $k$th row of ${\mathbf A}_{\tau(k)}$~. For a power index $(\alpha_1, \ldots, \alpha_m)$ with $\alpha_i \in \{0, 1, \ldots, d+J-1\}$ and $\sum_{i=1}^m \alpha_i = d+J-1$, we denote
$$\tau\in (\alpha_1, \ldots, \alpha_m)$$
if $\alpha_i = \#\{j : \tau(j)=i\}$ for each $i=1, \ldots, m$. In terms of the construction of ${\mathbf A}_\tau$, it says that $\alpha_i$ rows of ${\mathbf A}_\tau$ are from the matrix ${\mathbf A}_i$~.

According to the Leibniz formula for the determinant,
$$|{\mathbf F}| = \left|\sum_{i=1}^m n_i {\mathbf A}_i\right| = \sum_{\sigma \in S_{d+J-1}} (-1)^{{\rm sgn}(\sigma)} \prod_{k=1}^{d+J-1} \sum_{i=1}^m n_i a^{(i)}_{k,\sigma(k)}$$
where $\sigma$ is a permutation of $\{1, 2, \ldots, d+J-1\}$, and ${\rm sgn}(\sigma)$ is the sign or signature of $\sigma$. Therefore,
\begin{eqnarray*}
c_{\alpha_1, \ldots, \alpha_m} &=& \sum_{\sigma \in S_{d+J-1}} (-1)^{{\rm sgn}(\sigma)} \sum_{\tau \in (\alpha_1, \ldots, \alpha_m)} \prod_{k=1}^{d+J-1} a^{(\tau(k))}_{k,\sigma(k)}\\
&=& \sum_{\tau \in (\alpha_1, \ldots, \alpha_m)} \sum_{\sigma \in S_{d+J-1}} (-1)^{{\rm sgn}(\sigma)} \prod_{k=1}^{d+J-1} a^{(\tau(k))}_{k,\sigma(k)}\\
&=& \sum_{\tau \in (\alpha_1, \ldots, \alpha_m)}|{\mathbf A}_\tau|
\end{eqnarray*}
\hfill{$\Box$}

\bigskip\noindent
{\bf Proof of Lemma~\ref{coelemma3}}\quad
To simplify the notations, we let $i_s=s+1$, $s=0, \ldots, d$.
That is, $\alpha_1=J-1$, $\alpha_2 = \cdots = \alpha_{d+1} = 1$. There are only two types of $\tau \in (\alpha_1, \ldots, \alpha_m)$, such that, $|{\mathbf A}_\tau|$ may not be $0$.

{\it $\tau$ of type I:} There exist $1\leq k\leq d$, $2\leq l\leq d+1$, and $1\leq q\leq J-1$, such that, $\tau(k)=1$ and $\tau(d+q)=l$. Following a similar procedure as in the proof of Lemma~\ref{coelemma2}, we obtain
$$|{\mathbf A}_{\tau}| = \prod_{i=2}^{d+1} e_i \cdot |{\mathbf A}_{13}|\cdot (-1)^d |{\mathbf X}_{\bf 1}[1, 2,\ldots, d+1]| \cdot  (-1)^{{\rm sgn}(\tau)}\prod_{s=1}^d x_{\tau(s)s}\cdot \frac{c_{lq}}{e_l}$$

{\it $\tau$ of type II:} $\tau(d+1)=\cdots = \tau(d+J-1)=1$ and $\{\tau(1), \ldots,$ $\tau(d)\}$ $=$ $\{2, \ldots, d+1\}$. It can be verified that
$$|{\mathbf A}_{\tau}| = \prod_{i=2}^{d+1} e_i \cdot |{\mathbf A}_{13}|\cdot (-1)^d |{\mathbf X}_{\bf 1}[1, 2,\ldots, d+1]| \cdot  (-1)^{{\rm sgn}(\tau)}\prod_{s=1}^d x_{\tau(s)s}$$

According to Theorem~\ref{polynomialtheorem},
\begin{eqnarray*}
& & c_{\alpha_1, \ldots, \alpha_m} = \sum_{\tau \mbox{ of type I}} |{\mathbf A}_\tau| + \sum_{\tau \mbox{ of type II}} |{\mathbf A}_\tau|\\
&=& \prod_{i=2}^{d+1} e_i \cdot |{\mathbf A}_{13}|\cdot (-1)^d |{\mathbf X}_{\bf 1}[1, 2,\ldots, d+1]| \cdot
\left(\sum_{k=1}^d \sum_{l=2}^{d+1}\sum_{\tau \in S_{d+1} : \tau(k)=1, \tau(d+1)=l}\right.
\end{eqnarray*}
\begin{eqnarray*}
& &\left. (-1)^{{\rm sgn}(\tau)}\prod_{s=1}^d x_{\tau(s)s} \sum_{q=1}^{J-1} \frac{c_{lq}}{e_l}
+\sum_{\tau\in S_{d+1} : \tau(d+1)=1} (-1)^{{\rm sgn}(\tau)}\prod_{s=1}^d x_{\tau(s)s}\right)\\
&=& \prod_{i=2}^{d+1} e_i \cdot |{\mathbf A}_{13}|\cdot (-1)^d |{\mathbf X}_{\bf 1}[1, 2,\ldots, d+1]| \cdot
\sum_{\tau\in S_{d+1}} (-1)^{{\rm sgn}(\tau)}\prod_{s=1}^d x_{\tau(s)s}\\
&=& \prod_{i=2}^{d+1} e_i \cdot |{\mathbf A}_{13}|\cdot (-1)^d |{\mathbf X}_{\bf 1}[1, 2,\ldots, d+1]| \cdot (-1)^d |{\mathbf X}_{\bf 1}[1, 2,\ldots, d+1]|\\
&=& \prod_{i=2}^{d+1} e_i \cdot |{\mathbf A}_{13}|\cdot |{\mathbf X}_{\bf 1}[1, 2,\ldots, d+1]|^2
\end{eqnarray*}
where $S_{d+1}$ is the set of permutations of $\{1, \ldots, d+1\}$. The general case with $i_0, i_1, \ldots, i_d$ can be obtained similarly.
\hfill{$\Box$}

\bigskip\noindent
{\bf Proof of Theorem~\ref{mindesigntheorem2}}\quad
Suppose ${\rm Rank}({\mathbf X}_{\bf 1})=d+1$. Then there exist $i_0, \ldots,$ $i_d$ $\in$ $\{1, \ldots, m\}$, such that,
$|{\mathbf X}_{\bf 1}[i_0,$ $i_1,$ $\ldots,$ $i_d]|\neq 0$. According to Lemma~\ref{coelemma1}, $f({\bf p})$ can be regarded as an order-$(J-1)$ polynomial of $p_{i_0}$~. Let $p_{i_0}=x \in (0,1)$ and $p_i = (1-x)/(m-1)$ for $i\neq i_0$~. Based on Lemma~\ref{coelemma3}, $f({\bf p})$ can be written as
\begin{eqnarray*}
f_{i_0}(x) &=& a_{J-1} x^{J-1}\left(\frac{1-x}{m-1}\right)^d + a_{J-2} x^{J-2} \left(\frac{1-x}{m-1}\right)^{d+1}\\
& & + \cdots + a_1 x\left(\frac{1-x}{m-1}\right)^{d+J-2} + a_0 \left(\frac{1-x}{m-1}\right)^{d+J-1}, \mbox{  where}\\
a_{J-1} &=& |{\mathbf A}_{i_03}| \sum_{\{i_1', \ldots, i_d'\} \subset \{1, \ldots, m\}\setminus\{i_0\}}\> \prod_{s=1}^d e_{i_s'}\> |{\mathbf X}_{\bf 1}[i_0, i_1', \ldots, i_d']|^2 > 0
\end{eqnarray*}
Therefore, $\lim_{x\rightarrow 1^-} (1-x)^{-d} x^{1-J} f_{i_0}(x) = (m-1)^{-d} a_{J-1} > 0$. That is, $f({\bf p}) > 0$ for $p_{i_0}=x$ close enough to $1$ and $p_i=(1-x)/(m-1)$ for $i\neq i_0$~.

In order to justify that the condition ${\rm Rank}({\mathbf X}_{\bf 1})=d+1$ is also necessary, we only need to show that $f({\mathbf p})\equiv 0$ if ${\rm Rank}({\mathbf X}_{\bf 1}) \leq d$. Actually, for any $\tau : \{1, \ldots, d+J-1\} \rightarrow \{1, \ldots, m\}$, we construct ${\mathbf A}^{(1)}_\tau$ as in the proof of Lemma~\ref{coelemma2}. Then $|{\mathbf A}_\tau|=\prod_{s=1}^dx_{\tau(s)s}\cdot |{\mathbf A}^{(1)}_\tau|$. Similar as in the proof of Lemma~\ref{coelemma2}, for $t=1, \ldots, d$, we add $x_{\tau(1)t}(-e_{\tau(1)}, \ldots, -e_{\tau(d)}, c_{\tau(d+1)1}, $ $\ldots,$ $c_{\tau(d+J-1),J-1})^T$ to the $t$th column of ${\mathbf A}^{(1)}_\tau$. We denote the resulting matrix by ${\mathbf A}^{(3)}_\tau$. Note that $|{\mathbf A}^{(1)}_\tau| = |{\mathbf A}^{(3)}_\tau|$. We consider the sub-matrix ${\mathbf A}^{(3)}_{\tau d}$ which consists of the first $d$ columns of ${\mathbf A}^{(3)}_\tau$. For $s \in \tau^{-1}(\tau(1))$, the $s$th row of ${\mathbf A}^{(3)}_{\tau d}$ is simply $0$. For $s=2, \ldots, k$, the $s$th row of ${\mathbf A}^{(3)}_{\tau d}$ is $e_{\tau(s)}(x_{\tau(s)1}-x_{\tau(1)1}, \ldots, x_{\tau(s)d} - x_{\tau(1)d})$. For $s=1, \ldots, J-1$, the $(d+s)$th row of ${\mathbf A}^{(3)}_{\tau d}$ is $-c_{\tau(d+s)s}(x_{\tau(d+s)1}-x_{\tau(1)1}, \ldots, x_{\tau(d+s)d} - x_{\tau(1)d})$.
We claim that ${\rm Rank}({\mathbf A}^{(3)}_{\tau d}) \leq d-1$. Otherwise, if ${\rm Rank}({\mathbf A}^{(3)}_{\tau d})=d$, then there exist $i_1, \ldots, i_d \in \{2, \ldots, d+J-1\}$, such that, the sub-matrix consisting of the $i_1{\rm th}, \ldots, i_d$th rows of ${\mathbf A}^{(3)}_{\tau d}$ is nonsingular. Then the sub-matrix consisting of the $\tau(1){\rm th}, \tau(i_1){\rm th}, \ldots, \tau(i_d)$th rows of ${\mathbf X}_{\bf 1}$ is nonsingular, which implies ${\rm Rank}({\mathbf X}_{\bf 1}) = d+1$. The contradiction implies ${\rm Rank}({\mathbf A}^{(3)}_{\tau d}) \leq d-1$. Then $|{\mathbf A}^{(3)}_\tau|=0$ and thus $|{\mathbf A}_\tau|=0$ for each $\tau$. Based on Theorem~\ref{polynomialtheorem}, $|{\mathbf F}|\equiv 0$ and thus $f({\mathbf p})\equiv 0$.
\hfill{$\Box$}

\bigskip\noindent
{\bf Proof of Theorem~\ref{s+theorem}}\quad
Combining Theorem~\ref{sumFtheorem} and Theorem~\ref{mindesigntheorem2}, it is straightforward that $f({\mathbf p}) = 0$ if ${\rm Rank}($ ${\mathbf X}_{\bf 1}[\{i\mid p_i>0\}]) \leq d$. We only need to show that $f({\mathbf p}) > 0$ if ${\rm Rank}({\mathbf X}_{\bf 1}[\{i\mid p_i>0\}])=d+1$. Due to Theorem~\ref{sumFtheorem}, we only need to verify the case $p_i>0, i=1, \ldots, m$, because otherwise we may simply remove all support points with $p_i=0$.

Suppose $p_i > 0, i=1, \ldots, m$ and ${\rm Rank}({\mathbf X}_{\bf 1}) = d+1$. Then there exist $i_0, \ldots, i_d \in \{1, \ldots, m\}$, such that, $|{\mathbf X}_{\bf 1}[i_0, \ldots, i_d]|\neq 0$. According to the proof of Theorem~\ref{mindesigntheorem2}, for each $i \in \{i_0, \ldots, i_d\}$, there exists an $\epsilon_i \in (0,1)$, such that, $f({\mathbf p}) > 0$ as long as $p_i=x \in (1-\epsilon_i, 1)$ and $p_j = (1-x)/(m-1)$ for $j\neq i$. On the other hand, for each $i \notin \{i_0, \ldots, i_d\}$, if we denote the $j$th row of ${\mathbf X}_{\bf 1}$ by $\alpha_j$, $j=1, \ldots, m$, then $\alpha_i = a_0\alpha_{i_0} + \cdots + a_d \alpha_{i_d}$ for some real numbers $a_0, \ldots, a_d$~. Since $\alpha_i\neq 0$, then at least one $a_i\neq 0$. Without any loss of generality, we assume $a_0\neq 0$. Then it can be verified that $|{\mathbf X}_{\bf 1}[i, i_1, \ldots, i_d]|\neq 0$ too. Following the proof of Theorem~\ref{mindesigntheorem2} again, for such an $i \notin \{i_0, \ldots, i_d\}$, there also exists an $\epsilon_i \in (0,1)$, such that, $f({\mathbf p}) > 0$ as long as $p_i=x \in (1-\epsilon_i, 1)$ and $p_j = (1-x)/(m-1)$ for $j\neq i$. Let $\epsilon_* =\min\{\min_i \epsilon_i, (m-1)\min_i p_i, 1-1/m\}/2$. For $i=1, \ldots, m$, denote $\delta_i = (\delta_{i1}, \ldots, \delta_{im})^T \in S$ with $\delta_{ii}=1-\epsilon_*$ and $\delta_{ij}=\epsilon_*/(m-1)$ for $j\neq i$.
It can be verified that ${\mathbf p}=a_1 \delta_1 + \cdots + a_m \delta_m$ with $a_i = (p_i - \epsilon_*/(m-1))/(1-m\epsilon_*/(m-1))$. By the choice of $\epsilon_*$, $f(\delta_i)>0$, $a_i>0, i=1, \ldots, m$, and $\sum_i a_i = 1$. Then $f({\mathbf p})>0$ according to Lemma~\ref{|F|lemma}.
\hfill{$\Box$}

\bigskip\noindent
{\bf Proof of Corollary~\ref{d=1J=3m=3corollary}}\quad
In order to check when a minimally supported design supported only on $\{x_1, x_2\}$ is D-optimal, we add one more support point, that is, $x_3$~. According to Theorem~\ref{polynomialtheorem}, Lemmas~\ref{coelemma1}, \ref{coelemma2}, and \ref{coelemma3}, the objective function for a D-optimal approximate design on $\{x_1, x_2, x_3\}$ is
$f(p_1, p_2, p_3)=p_1 p_2 (c_{210} p_1 + c_{120} p_2) + p_1 p_3 (c_{201} p_1 + c_{102} p_3) + p_2 p_3 (c_{021} p_2 + c_{012} p_3) + c_{111} p_1 p_2 p_3$,
where
\begin{eqnarray*}
c_{210} &=& e_2 g_{11}^2 g_{12}^2 (\pi_{11}\pi_{12}\pi_{13})^{-1} (x_1 - x_2)^2 > 0\\
c_{120} &=& e_1 g_{21}^2 g_{22}^2 (\pi_{21}\pi_{22}\pi_{23})^{-1} (x_1 - x_2)^2 > 0\\
c_{201} &=& e_3 g_{11}^2 g_{12}^2 (\pi_{11}\pi_{12}\pi_{13})^{-1} (x_1 - x_3)^2 > 0\\
c_{102} &=& e_1 g_{31}^2 g_{32}^2 (\pi_{31}\pi_{32}\pi_{33})^{-1} (x_1 - x_3)^2 > 0\\
c_{021} &=& e_3 g_{21}^2 g_{22}^2 (\pi_{21}\pi_{22}\pi_{23})^{-1} (x_2 - x_3)^2 > 0\\
c_{012} &=& e_2 g_{31}^2 g_{32}^2 (\pi_{31}\pi_{32}\pi_{33})^{-1} (x_2 - x_3)^2 > 0\\
c_{111} &=& e_1 (u_{22}u_{31} + u_{21}u_{32} - 2b_{22}b_{32})(x_1 - x_2)(x_1 - x_3) + \\
& & e_2 (u_{12}u_{31} + u_{11} u_{32} - 2 b_{12} b_{32}) (x_2 - x_1) (x_2 - x_3) + \\
& & e_3 (u_{12} u_{21} + u_{11}u_{22} - 2b_{12} b_{22}) (x_3-x_1)(x_3-x_2)
\end{eqnarray*}
Based on Theorem~\ref{KKTtheorem}, the design ${\mathbf p} = (p_1^*, p_2^*, 0)^T$ is D-optimal if and only if $$\partial f({\mathbf p})/\partial f(p_1) = \partial f({\mathbf p})/\partial f(p_2) \geq \partial f({\mathbf p})/\partial f(p_3)$$
Similar conclusions could be justified for $x_4, \ldots, x_m$ if $m\geq 4$.
\hfill{$\Box$}

\bigskip\noindent
{\bf Proof of Theorem~\ref{f4theorem}}\quad
According to the solutions provided by the software {\tt Mathematica}, the largest root of equation~(\ref{y1equation}) after simplification is
\begin{equation}\label{y1solution}
y_1=-\frac{b_2}{3} - \frac{2^{1/3} (3 b_1 - b_2^2)}{3 A^{1/3}} +  \frac{A^{1/3}}{3\times 2^{1/3}}
\end{equation}
where $A=-27 b_0 + 9 b_1 b_2 - 2 b_2^3 + 3^{3/2} (27 b_0^2 + 4 b_1^3 - 18 b_0 b_1 b_2 - b_1^2 b_2^2 + 4 b_0 b_2^3)^{1/2}$, and $b_i = c_i/c_3$, $i=0,1,2$.
Note that the calculation of $A$ and thus $y_1$ should be regarded as operations among complex numbers since the expression under square root could be negative. Nevertheless, $y_1$ at the end would be a real number. Thus we are able to provide the analytic solution maximizing $f(p_1, p_2, p_3)$.
\hfill{$\Box$}

\bigskip\noindent
{\bf Proof of Corollary~\ref{d=2J=3m=4corollary}}\quad
In order to check when a minimally supported design is D-optimal, we first add the four design points, that is, we consider four design points $(x_{i1}, x_{i2})$, $i=1,2,3,4$ and check when the D-optimal design could be constructed on the first three design points. Let ${\mathbf X}_{\bf 1}$ be defined as in Lemma~\ref{coelemma3}. In this case, ${\mathbf X}_{\bf 1}$ is a $4\times 3$ matrix.
Following Theorem~\ref{polynomialtheorem}, Lemmas~\ref{coelemma1}, \ref{coelemma2}, and \ref{coelemma3}, the objective function for a minimally supported design at $(d, J, m) = (2, 3, 4)$ is
\begin{eqnarray*}
f(p_1, p_2, p_3, p_4) &=& c_{1111} p_1 p_2 p_3 p_4 \\
&+& |{\mathbf X}_{\bf 1}[1,2,3]|^2 e_1 e_2 e_3 \cdot p_1 p_2 p_3 (w_1 p_1 + w_2 p_2 + w_3 p_3)\\
&+& |{\mathbf X}_{\bf 1}[1,2,4]|^2 e_1 e_2 e_4 \cdot p_1 p_2 p_4 (w_1 p_1 + w_2 p_2 + w_4 p_4)\\
&+& |{\mathbf X}_{\bf 1}[1,3,4]|^2 e_1 e_3 e_4 \cdot p_1 p_3 p_4 (w_1 p_1 + w_3 p_3 + w_4 p_4)\\
&+& |{\mathbf X}_{\bf 1}[2,3,4]|^2 e_2 e_3 e_4 \cdot p_2 p_3 p_4 (w_2 p_2 + w_3 p_3 + w_4 p_4)
\end{eqnarray*}
where $e_i = u_{i1} + u_{i2} - 2b_{i2}$, $w_i = e_i^{-1} g_{i1}^2 g_{i2}^2 (\pi_{i1}\pi_{i2}\pi_{i3})^{-1}$, $i=1, 2, 3, 4$, and
\begin{equation}\label{c1111}
c_{1111} = \sum_{1\leq i < j \leq 4} e_i e_j (u_{k1} u_{l2} + u_{k2} u_{l1} - 2 b_{k2} b_{l2})\cdot |{\mathbf X}_{\bf 1}[i,j,k]|\cdot |{\mathbf X}_{\bf 1}[i,j,l]|
\end{equation}
with $\{i,j,k,l\}=\{1,2,3,4\}$ given $1\leq i<j\leq 4$.

According to Theorem~\ref{KKTtheorem}, a minimally supported design ${\mathbf p} = (p_1^*, p_2^*, p_3^*,$ $0)^T$ in this case is D-optimal if and only if
$\partial f/\partial p_1 = \partial f/\partial p_2  = \partial f/\partial p_3 \geq \partial f/\partial p_4$
at ${\mathbf p}$. Then $\partial f/\partial p_1 = \partial f/\partial p_2  = \partial f/\partial p_3$ at ${\mathbf p}$ is equivalent to (1) of Corollary~\ref{d=2J=3m=4corollary}, and
$\partial f/\partial p_4 \leq \partial f/\partial p_1$ at ${\mathbf p}$ leads to (2) of Corollary~\ref{d=2J=3m=4corollary} since the forms of $\partial f/\partial p_i$ at ${\mathbf p}$, $i=1,2,3$
will not change if more than four design points (i.e., $m > 4$) are added into consideration. Note that\\
$|{\mathbf X}_{\bf 1}[1,2,3]|^2 e_1 e_2 e_3 p_2^* p_3^* (2w_1 p_1^* + w_2 p_2^* + w_3 p_3^*)$ in (2) of Corollary~\ref{d=2J=3m=4corollary} is equal to $\partial f/\partial p_1$ at ${\mathbf p}$. It could be replaced with $|{\mathbf X}_{\bf 1}[1,2,3]|^2 e_1 e_2 e_3 p_1^* p_3^* (w_1 p_1^* + 2w_2 p_2^* + w_3 p_3^*)$ (i.e., $\partial f/\partial p_2$), or $|{\mathbf X}_{\bf 1}[1,2,3]|^2 e_1 e_2 e_3 p_1^* p_2^* (w_1 p_1^* + w_2 p_2^* + 2 w_3 p_3^*)$ (i.e., $\partial f/\partial p_3$), since these three are all equal.
\hfill{$\Box$}

\subsection{Maximization of $f_i(z)$ in Section 3}\label{S4section}

According to Theorem~\ref{f_i(z)theorem}, $f_i(z)$ is an order-$(d+J-1)$ polynomial of $z$. In other to determine its coefficients $a_0, a_1, \ldots, a_{J-1}$ as in (\ref{f_i(z)poly}), we need to calculate $f_i(0),$ $f_i(1/2),$ $f_i(1/3),$ $\ldots,$ $f_i(1/J)$, which are $J$ determinants defined in (\ref{f_i(z)}).

Note that ${\mathbf B}_{J-1}^{-1}$ is a matrix determined by $J-1$ only. For example, $B_1^{-1} = 1$ for $J=2$,
\[
B_2^{-1}=\left(
\begin{array}{rr}
2 & -1\\
-1 & 1
\end{array}
\right),
B_3^{-1} = \left(
\begin{array}{rrr}
3 & -3 & 1\\
-\frac{5}{2} & 4 & -\frac{3}{2}\\
\frac{1}{2} & -1 & \frac{1}{2}
\end{array}
\right),
\]
\[
B_4^{-1} = \left(
\begin{array}{rrrr}
4 & -6 & 4 & -1\\
-\frac{13}{3} & \frac{19}{2} & -7 & \frac{11}{6}\\
\frac{3}{2} & -4 & \frac{7}{2} & -1\\
-\frac{1}{6} & \frac{1}{2} & -\frac{1}{2} & \frac{1}{6}
\end{array}
\right)
\]
for $J=3$, $4$, or $5$ respectively.

Once $a_0, \ldots, a_{J-1}$ in (\ref{f_i(z)poly}) are determined, the maximization of $f_i(z)$ on $z\in [0, 1]$ is numerically straightforward since it is a polynomial and its derivative $f_i'(z)$ is given by
\begin{equation}\label{f_i(z)'}
(1-z)^d \sum_{j=1}^{J-1} j a_j z^{j-1} (1-z)^{J-1-j} - (1-z)^{d-1} \sum_{j=0}^{J-1} (d+J-1-j) a_j z^j (1-z)^{J-1-j}
\end{equation}

\subsection{Exchange algorithm for D-optimal exact allocation in Section 4}\label{S5section}

\bigskip\noindent
{\it Exchange algorithm for D-optimal allocation $(n_1, \ldots, n_m)^T$ given $n>0$:}
\begin{itemize}
 \item[$1^\circ$] Start with an initial design ${\mathbf n}=(n_1,\ldots,n_m)^T$ such that $f({\mathbf n})>0$.
 \item[$2^\circ$] Set up a random order of $(i,j)$ going through all pairs $\{(1,2),$ $(1,3),$ $\ldots,$ $(1,m),$ $(2,3),$ $\ldots,(m-1, m)\}$.
 \item[$3^\circ$] For each $(i,j)$, let $c=n_i+n_j$~. If $c=0$, let ${\mathbf n}^*_{ij}={\mathbf n}$. Otherwise, there are two cases. {\it Case one:} $0<c\leq J$, we calculate $f_{ij}(z)$ as defined in (\ref{f_{ij}}) for $z=0, 1, \ldots, c$ directly and find $z^*$ which maximizes $f_{ij}(z)$. {\it Case two:} $c > J$, we first calculate $f_{ij}(z)$ for $z=0, 1, \ldots, J$; secondly determine $c_0, c_1, \ldots, c_J$ in (\ref{f_ij(z)poly}) according to Theorem~\ref{f_ij(z)theorem}; thirdly calculate $f_{ij}(z)$ for $z=J+1, \ldots, c$ based on (\ref{f_ij(z)poly}); fourthly find $z^*$ maximizing $f_{ij}(z)$ for $z=0, \ldots, c$. For both cases, we define
     $${\mathbf n}^*_{ij} =\left(n_1,\ldots,n_{i-1},z^*,n_{i+1},\ldots,n_{j-1},c-z^*,n_{j+1},\ldots,n_m\right)^T$$
     Note that $f({\mathbf n}^*_{ij})=f_{ij}(z^*)\geq f({\mathbf n})>0$. If $f({\mathbf n}^*_{ij}) > f({\mathbf n})$, replace ${\mathbf n}$ with ${\mathbf n}^*_{ij}$~, and $f({\mathbf n})$ with $f({\mathbf n}^*_{ij})$.
 \item[$4^\circ$] Repeat $2^\circ\sim 3^\circ$ until convergence, that is, $f({\mathbf n}^*_{ij})= f({\mathbf n})$ in step~$3^\circ$ for any $(i,j)$.
\end{itemize}

\subsection{More on Example~\ref{industrialexample}: Polysilicon Deposition Study}\label{S6section}

Table~\ref{polysilicontable} shows the list of experimental settings for the polysilicon deposition study. The factors are {\tt decomposition temperature}($A$), {\tt decomposition pressure}($B$),  {\tt nitrogen} {\tt flow} ($C$), {\tt silane flow}($D$), {\tt setting time}($E$), {\tt cleaning method}($F$). Column~1 provides original indices of experimental settings out of 729 distinct ones. For each experimental setting labelled ``1" in a design, 9 responses are collected (\cite{phadke1989}) and assumed to be independent.

\begin{table}[h!]\caption{Polysilicon Deposition Study: Experimental Settings for the Original, Rounded Approximate, and D-optimal Exact Designs}
\label{polysilicontable}
\vskip -0.3in
\begin{center}
\tiny
\begin{tabular}{|r|rrrrrr|c|c|c|}  \hline
Index & $A$ & $B$ & $C$ & $D$ & $E$ & $F$ & Original & Rounded & D-optimal\\ \hline
1   &  1   &  1   &  1   &  1   &  1   &  1   &  1   &  0   &  0 \\
76   &  1   &  1   &  3   &  3   &  2   &  1   &  1   &  0   &  0 \\
89   &  1   &  2   &  1   &  1   &  3   &  2   &  1   &  0   &  0 \\
98   &  1   &  2   &  1   &  2   &  3   &  2   &  0   &  0   &  1 \\
111   &  1   &  2   &  2   &  1   &  1   &  3   &  0   &  0   &  1 \\
116   &  1   &  2   &  2   &  1   &  3   &  2   &  0   &  1   &  0 \\
122   &  1   &  2   &  2   &  2   &  2   &  2   &  1   &  0   &  0 \\
130   &  1   &  2   &  2   &  3   &  2   &  1   &  0   &  0   &  1 \\
167   &  1   &  3   &  1   &  1   &  2   &  2   &  0   &  0   &  1 \\
181   &  1   &  3   &  1   &  3   &  1   &  1   &  0   &  1   &  0 \\
199   &  1   &  3   &  2   &  2   &  1   &  1   &  0   &  1   &  1 \\
201   &  1   &  3   &  2   &  2   &  1   &  3   &  1   &  0   &  0 \\
243   &  1   &  3   &  3   &  3   &  3   &  3   &  1   &  0   &  1 \\
258   &  2   &  1   &  1   &  2   &  2   &  3   &  1   &  0   &  0 \\
286   &  2   &  1   &  2   &  2   &  3   &  1   &  0   &  1   &  0 \\
290   &  2   &  1   &  2   &  3   &  1   &  2   &  1   &  0   &  0 \\
291   &  2   &  1   &  2   &  3   &  1   &  3   &  0   &  1   &  0 \\
294   &  2   &  1   &  2   &  3   &  2   &  3   &  0   &  0   &  1 \\
299   &  2   &  1   &  3   &  1   &  1   &  2   &  0   &  0   &  1 \\
301   &  2   &  1   &  3   &  1   &  2   &  1   &  0   &  1   &  0 \\
313   &  2   &  1   &  3   &  2   &  3   &  1   &  0   &  0   &  1 \\
331   &  2   &  2   &  1   &  1   &  3   &  1   &  0   &  1   &  1 \\
336   &  2   &  2   &  1   &  2   &  1   &  3   &  0   &  1   &  1 \\
339   &  2   &  2   &  1   &  2   &  2   &  3   &  0   &  1   &  0 \\
350   &  2   &  2   &  1   &  3   &  3   &  2   &  0   &  1   &  0 \\
365   &  2   &  2   &  2   &  2   &  2   &  2   &  0   &  0   &  1 \\
376   &  2   &  2   &  2   &  3   &  3   &  1   &  1   &  0   &  0 \\
384   &  2   &  2   &  3   &  1   &  2   &  3   &  1   &  0   &  0 \\
394   &  2   &  2   &  3   &  2   &  3   &  1   &  0   &  1   &  0 \\
399   &  2   &  2   &  3   &  3   &  1   &  3   &  0   &  1   &  0 \\
407   &  2   &  3   &  1   &  1   &  1   &  2   &  0   &  0   &  1 \\
421   &  2   &  3   &  1   &  2   &  3   &  1   &  1   &  0   &  0 \\
461   &  2   &  3   &  3   &  1   &  1   &  2   &  1   &  1   &  0 \\
464   &  2   &  3   &  3   &  1   &  2   &  2   &  0   &  1   &  0 \\
495   &  3   &  1   &  1   &  1   &  3   &  3   &  0   &  1   &  0 \\
501   &  3   &  1   &  1   &  2   &  2   &  3   &  0   &  0   &  1 \\
505   &  3   &  1   &  1   &  3   &  1   &  1   &  0   &  0   &  1 \\
521   &  3   &  1   &  2   &  1   &  3   &  2   &  0   &  0   &  1 \\
522   &  3   &  1   &  2   &  1   &  3   &  3   &  1   &  0   &  0 \\
536   &  3   &  1   &  2   &  3   &  2   &  2   &  0   &  1   &  0 \\
557   &  3   &  1   &  3   &  2   &  3   &  2   &  1   &  0   &  0 \\
558   &  3   &  1   &  3   &  2   &  3   &  3   &  0   &  1   &  0 \\
569   &  3   &  2   &  1   &  1   &  1   &  2   &  0   &  1   &  0 \\
588   &  3   &  2   &  1   &  3   &  1   &  3   &  1   &  0   &  0 \\
625   &  3   &  2   &  3   &  1   &  2   &  1   &  0   &  0   &  1 \\
631   &  3   &  2   &  3   &  2   &  1   &  1   &  1   &  0   &  0 \\
641   &  3   &  2   &  3   &  3   &  1   &  2   &  0   &  0   &  1 \\
671   &  3   &  3   &  1   &  3   &  2   &  2   &  1   &  0   &  0 \\
679   &  3   &  3   &  2   &  1   &  2   &  1   &  1   &  0   &  0 \\
\hline
\end{tabular}
\normalsize
\end{center}
\end{table}

Table~\ref{polysilicontable2} shows the model matrix for the D-optimal design ${\mathbf n}_o$ found for the polysilicon deposition study. In this table, each 3-level factor is represented by its linear component and quadratic component. Thus there are level combinations of 12 predictors.

\begin{table}[h!]\caption{Polysilicon Deposition Study: Model Matrix for the D-optimal Design}
\label{polysilicontable2}
\begin{center}
\begin{tabular}{|r|rrrrrrrrrrrr|}  \hline
Index  & $A_1$ & $A_2$ & $B_1$ & $B_2$ & $C_1$ & $C_2$ & $D_1$ & $D_2$ & $E_1$ & $E_2$ & $F_1$ & $F_2$ \\ \hline
 98   &  $-1$   &  $ 1$   &  $ 0$   &  $-2$   &  $-1$   &  $ 1$   &  $ 0$   &  $-2$   &  $ 1$   &  $ 1$   &  $ 0$   &  $-2$ \\
111   &  $-1$   &  $ 1$   &  $ 0$   &  $-2$   &  $ 0$   &  $-2$   &  $-1$   &  $ 1$   &  $-1$   &  $ 1$   &  $ 1$   &  $ 1$ \\
130   &  $-1$   &  $ 1$   &  $ 0$   &  $-2$   &  $ 0$   &  $-2$   &  $ 1$   &  $ 1$   &  $ 0$   &  $-2$   &  $-1$   &  $ 1$ \\
167   &  $-1$   &  $ 1$   &  $ 1$   &  $ 1$   &  $-1$   &  $ 1$   &  $-1$   &  $ 1$   &  $ 0$   &  $-2$   &  $ 0$   &  $-2$ \\
199   &  $-1$   &  $ 1$   &  $ 1$   &  $ 1$   &  $ 0$   &  $-2$   &  $ 0$   &  $-2$   &  $-1$   &  $ 1$   &  $-1$   &  $ 1$ \\
243   &  $-1$   &  $ 1$   &  $ 1$   &  $ 1$   &  $ 1$   &  $ 1$   &  $ 1$   &  $ 1$   &  $ 1$   &  $ 1$   &  $ 1$   &  $ 1$ \\
294   &  $ 0$   &  $-2$   &  $-1$   &  $ 1$   &  $ 0$   &  $-2$   &  $ 1$   &  $ 1$   &  $ 0$   &  $-2$   &  $ 1$   &  $ 1$ \\
299   &  $ 0$   &  $-2$   &  $-1$   &  $ 1$   &  $ 1$   &  $ 1$   &  $-1$   &  $ 1$   &  $-1$   &  $ 1$   &  $ 0$   &  $-2$ \\
313   &  $ 0$   &  $-2$   &  $-1$   &  $ 1$   &  $ 1$   &  $ 1$   &  $ 0$   &  $-2$   &  $ 1$   &  $ 1$   &  $-1$   &  $ 1$ \\
331   &  $ 0$   &  $-2$   &  $ 0$   &  $-2$   &  $-1$   &  $ 1$   &  $-1$   &  $ 1$   &  $ 1$   &  $ 1$   &  $-1$   &  $ 1$ \\
336   &  $ 0$   &  $-2$   &  $ 0$   &  $-2$   &  $-1$   &  $ 1$   &  $ 0$   &  $-2$   &  $-1$   &  $ 1$   &  $ 1$   &  $ 1$ \\
365   &  $ 0$   &  $-2$   &  $ 0$   &  $-2$   &  $ 0$   &  $-2$   &  $ 0$   &  $-2$   &  $ 0$   &  $-2$   &  $ 0$   &  $-2$ \\
407   &  $ 0$   &  $-2$   &  $ 1$   &  $ 1$   &  $-1$   &  $ 1$   &  $-1$   &  $ 1$   &  $-1$   &  $ 1$   &  $ 0$   &  $-2$ \\
501   &  $ 1$   &  $ 1$   &  $-1$   &  $ 1$   &  $-1$   &  $ 1$   &  $ 0$   &  $-2$   &  $ 0$   &  $-2$   &  $ 1$   &  $ 1$ \\
505   &  $ 1$   &  $ 1$   &  $-1$   &  $ 1$   &  $-1$   &  $ 1$   &  $ 1$   &  $ 1$   &  $-1$   &  $ 1$   &  $-1$   &  $ 1$ \\
521   &  $ 1$   &  $ 1$   &  $-1$   &  $ 1$   &  $ 0$   &  $-2$   &  $-1$   &  $ 1$   &  $ 1$   &  $ 1$   &  $ 0$   &  $-2$ \\
625   &  $ 1$   &  $ 1$   &  $ 0$   &  $-2$   &  $ 1$   &  $ 1$   &  $-1$   &  $ 1$   &  $ 0$   &  $-2$   &  $-1$   &  $ 1$ \\
641   &  $ 1$   &  $ 1$   &  $ 0$   &  $-2$   &  $ 1$   &  $ 1$   &  $ 1$   &  $ 1$   &  $-1$   &  $ 1$   &  $ 0$   &  $-2$ \\
\hline
\end{tabular}
\normalsize
\end{center}
\end{table}


{\large \bf References}
\begin{thebibliography}{11}
\expandafter\ifx\csname
natexlab\endcsname\relax\def\natexlab#1{#1}\fi
\expandafter\ifx\csname url\endcsname\relax
  \def\url#1{\texttt{#1}}\fi
\expandafter\ifx\csname urlprefix\endcsname\relax\def\urlprefix{URL
}\fi

\markboth{\hfill{\footnotesize\rm JIE YANG, LIPING TONG AND ABHYUDAY MANDAL} \hfill}
{\hfill {\footnotesize\rm D-OPTIMAL DESIGNS WITH ORDERED CATEGORICAL DATA} \hfill}

\bibitem[Abebe et~al.(2014)]{abebe2014}
Abebe, H.~T., Tan, F.~E., Van Breukelen, G.~J., Serroyen, J., and Berger, M.~P. (2014).
On the choice of a prior for Bayesian D-optimal designs for the logistic regression model with a single predictor. \textit{Communications in Statistics - Simulation and Computation}, \textbf{43}, 1811-1824.

\bibitem[Agresti(2013)]{agresti2013}
Agresti, A. (2013).
\textit{Categorical Data Analysis}, Third Edition. Wiley, New~Jersey.

\bibitem[Atkinson et~al.(2007)]{atkinson2007}
Atkinson, A.~C., Donev, A.~N.~and Tobias, R.~D. (2007).
\textit{Optimum Experimental Designs, with SAS}. Oxford University Press, New~York.

\bibitem[Chaloner and Larntz(1989)]{chaloner1989}
Chaloner, K. and Larntz, K. (1989).
Optimal Bayesian design applied to logistic regression experiments.
\textit{Journal of Statistical Planning and Inference}, \textbf{21}, 191-208.

\bibitem[Chaloner and Verdinelli(1995)]{chaloner1995}
Chaloner, K. and Verdinelli, I. (1995).
Bayesian experimental design: a review.
\textit{Statistical Science}, \textbf{10}, 273-304.

\bibitem[Chernoff(1953)]{chernoff1953}
Chernoff, H. (1953).
Locally optimal designs for estimating parameters.
\textit{Annals of Mathematical Statistics}, \textbf{24}, 586-602.

\bibitem[Christensen(2015)]{christensen2013}
Christensen, R.~H.~B. (2015).
Analysis of ordinal data with cumulative link models -- estimation with the R-package ordinal.
Available via \url{http://cran.r-project.org/web/packages/ordinal/vignettes/clm\_intro.pdf}

\bibitem[Dobson and Barnett(2008)]{dobson2008}
Dobson, A.~J. and Barnett, A. (2008).
\textit{An Introduction to Generalized Linear Models}, Third Edition. Chapman \& Hall/CRC, London.

\bibitem[Fedorov(1972)]{fedorov1972}
Fedorov, V.~V. (1972).
\textit{Theory of Optimal Experiments}. Academic Press, New~York.

\bibitem[Fedorov and Leonov(2014)]{fedorov2014}
Fedorov, V.~V. and Leonov, S.~L. (2014).
\textit{Optimal Design for Nonlinear Response Models}. Chapman \& Hall/CRC, New~York.


\bibitem[Ford et~al.(1989)]{ford1989}
Ford, I., Titterington, D.~M., and Kitsos, C.~P. (1989).
Recent advances in nonlinear experimental design.
\textit{Technometrics}, \textbf{31}, 49-60.

\bibitem[Ford, Torsney, and Wu(1992)]{ford1992}
Ford, I., Torsney, B., and Wu, C.~F.~J. (1992).
The use of a canonical form in the construction of locally optimal designs for non-linear problems.
\textit{Journal of the Royal Statistical Society, Series B}, \textbf{54}, 569-583.


\bibitem[Imhof(2001)]{imhof2001}
Imhof, L.~A. (2001).
Maximin designs for exponential growth models and heteroscedastic polynomial models.
\textit{Annals of Statistics}, \textbf{29}, 561-576.

\bibitem[Imhof, Lopez-Fidalgo, and Wong(2001)]{imhofwong2001}
Imhof, L., Lopez-Fidalgo, J., and Wong, W.K. (2001).
Efficiencies of rounded optimal approximate designs for small samples.
\textit{Statistica Neerlandica}, \textbf{55}, 301-318.

\bibitem[Joseph and Wu(2004)]{joseph2004}
Joseph, V.~R. and Wu, C.~F.~J. (2004).
Failure amplification method: an information maximization approach to categorical response optimization (with discussions).
\textit{Technometrics}, \textbf{46}, 1-31.


\bibitem[Karush(1939)]{karush1939}
Karush, W. (1939).
Minima of functions of several variables with inequalities as side constraints, M.Sc. Dissertation,
Department of Mathematics, University of Chicago.

\bibitem[Kiefer(1971)]{kiefer1971}
Kiefer, J. (1971).
The role of symmetry and approximation in exact design optimality.
In \textit{Statistical Decision Theory and Related Topics}, S.S. Gupta and J. Yackel (eds), 109-118, Academic Press, New~York.


\bibitem[Kiefer(1974)]{kiefer1974}
Kiefer, J. (1974).
General equivalence theory for optimum designs (approximate theory).
\textit{Annals of Statistics}, \textbf{2}, 849-879.

\bibitem[Khuri et~al.(2006)]{khuri2006}
Khuri, A.~I., Mukherjee, B., Sinha, B.~K. and Ghosh, M. (2006).
Design issues for generalized linear models: A review.
\textit{Statistical Science}, \textbf{21}, 376-399.

\bibitem[Krause et~al.(2001)]{krause2001}
Krause, M.~S., Madden, L.~V. and Hoitink, H.~A.~J. (2001).
Effect of potting mix microbial carrying capacity on biological control of Rhizoctonia damping-off of radish and Rhizoctonia crown and root rot of Poinsettia.
\textit{Phytopathology}, \textbf{91}, 1116-1123.

\bibitem[Kuhn and Tucker(1951)]{kuhn1951}
Kuhn, H.~W. and Tucker, A.~W. (1951).
Nonlinear programming.
\textit{Proceedings of 2nd Berkeley Symposium}, Berkeley: University of California Press, 481-492.

\bibitem[Liu and Agresti(2005)]{agresti2005}
Liu, I. and Agresti, A. (2005).
The analysis of ordered categorical data: An overview and a survey of recent developments.
\textit{Test}, \textbf{14}, 1-73.

\bibitem[McCullagh(1980)]{pmcc1980}
McCullagh, P. (1980).
Regression models for ordinal data.
\textit{Journal of the Royal Statistical Society, Series B}, \textbf{42}, 109-142.

\bibitem[McCullagh and Nelder(1989)]{pmcc1989}
McCullagh, P. and Nelder, J. (1989).
\textit{Generalized Linear Models}, Second Edition. Chapman and Hall/CRC, Boca~Raton.

\bibitem[Omer et~al.(2000)]{omer2000}
Omer, M.~A., Johnson, D.~A. and Rowe, R.~C. (2000).
Recovery of {\it Verticillium dahliae} from North American certified seed potatoes and
characterization of strains by vegetative compatibility and aggressiveness.
\textit{American Journal of Potato Research}, \textbf{77}, 325-331.

\bibitem[Osterstock et~al.(2010)]{osterstock2010}
Osterstock, J.~B., MacDonald, J.~C., Boggess, M.~M. and Brown, M.~S. (2010).
Analysis of ordinal outcomes from carcass data in beef cattle research.
\textit{Journal of Animal Science}, \textbf{88}, 3384-3389.

\bibitem[Perevozskaya et~al.(2003)]{perevozskaya2003}
Perevozskaya, I., Rosenberger, W.~F. and Haines, L.~M. (2003).
Optimal design for the proportional odds model.
\textit{The Canadian Journal of Statistics}, \textbf{31}, 225-235.

\bibitem[Phadke(1989)]{phadke1989}
Phadke, M.~S. (1989).
\textit{Quality Engineering using Robust Design}. Prentice-Hall, Englewood Cliffs.


\bibitem[Pronzato and Walter(1988)]{pronzato1988}
Pronzato, L. and Walter, E. (1988).
Robust experiment design via maximin optimization.
\textit{Mathematical Biosciences}, \textbf{89}, 161-176.


\bibitem[Pukelsheim(1993)]{pukelsheim1993}
Pukelsheim, F. (1993).
\textit{Optimal Design of Experiments}. John Wiley \& Sons, New~York.

\bibitem[Pukelsheim and Rieder(1992)]{pukelsheim1992}
Pukelsheim, F. and Rieder, S. (1992).
Efficient rounding of approximate designs.
\textit{Biometrika}, \textbf{79}, 763-770.

\bibitem[Randall(1989)]{randall1989}
Randall, J. (1989).
The analysis of sensory data by generalised linear model.
\textit{Biometrical Journal}, \textbf{31}, 781-793.

\bibitem[Silvey(1980)]{silvey1980}
Silvey, S.~D. (1980).
\textit{Optimal Design}. Chapman and Hall, London.

\bibitem[Song and Wong(1998)]{song1998}
Song, D. and Wong, W.~K. (1998).
Optimal two-point designs for the Michaelis-Menten model with heteroscedastic errors.
\textit{Communications in Statistics - Theory and Methods}, \textbf{27}, 1503-1516.

\bibitem[Stufken and Yang(2012)]{stufken2012}
Stufken, J. and Yang, M. (2012).
Optimal designs for generalized linear models.
In: \textit{Design and Analysis of Experiments}, Volume 3: Special Designs and Applications, K. Hinkelmann (ed.). Wiley, New~York.

\bibitem[Tong, Volkmer, and Yang(2014)]{tong2014}
Tong, L., Volkmer, H.~W., and Yang, J. (2014).
Analytic solutions for D-optimal factorial designs under generalized linear models.
\textit{Electronic Journal of Statistics}, \textbf{8}, 1322-1344.

\bibitem[Woods et~al.(2006)]{woods2006}
Woods, D.~C., Lewis, S.~M., Eccleston, J.~A. and Russell, K.~G. (2006).
Designs for generalized linear models with several variables and model uncertainty.
\textit{Technometrics}, \textbf{48}, 284-292.

\bibitem[Wu(2008)]{wu2008}
Wu, F-C. (2008).
Simultaneous optimization of robust design with quantitative and ordinal data.
\textit{International Journal of Industrial Engineering: Theory, Applications and Practice}, \textbf{15}, 231-238.


\bibitem[Wu and Hamada(2009)]{wu2009}
Wu, C.~F.~J. and Hamada, M. (2009).
\textit{Experiments: Planning, Analysis, and Optimization}, Second Edition. Wiley, New~York.

\bibitem[Yang and Mandal(2015)]{ym2014}
Yang, J. and Mandal, A. (2015).
D-optimal factorial designs under generalized linear models.
\textit{Communications in Statistics - Simulation and Computation}, \textbf{44}, 2264-2277.

\bibitem[Yang, Mandal, and Majumdar(2016)]{ymm2013}
Yang, J., Mandal, A., and Majumdar, D. (2016).
Optimal designs for $2^k$ factorial experiments with binary response.
\textit{Statistica Sinica}, \textbf{26}, 385-411.

\bibitem[Zocchi and Atkinson(1999)]{atkinson1999}
Zocchi S.~S. and Atkinson, A.~C. (1999).
Optimum experimental designs for multinomial logistic models.
\textit{Biometrics}, \textbf{55}, 437-444.

\end{thebibliography}
\end{document}